\documentclass[12pt]{amsart}

\usepackage{amsmath}
\usepackage[all]{xypic}
\usepackage[centering, hmargin=1in, vmargin=1in]{geometry}
\usepackage{amssymb}

\newtheorem{lemma}{Lemma}[section]
\newtheorem{theorem}[lemma]{Theorem}
\newtheorem{corollary}[lemma]{Corollary}
\theoremstyle{definition}
\newtheorem{definition}[lemma]{Definition}
\newtheorem{definitions}[lemma]{Definitions}

\newtheorem{remark}[lemma]{Remark}
\newtheorem{prop}[lemma]{Proposition}
\newtheorem{acknowledgements}{Acknowledgements}

\newcommand{\mr}{\mathrm}
\DeclareMathOperator{\kum}{Kum}

\author{Afsaneh Mehran}
\title{Double covers of Kummer surfaces}
\address{Section de Math\'ematiques, Universit\'e de Gen\`eve, 2-4 rue du Li\`evre, Case Postale 64, 1211 Gen\'eve 4, Suisse.\\}
\email{Afsaneh.Mehran@math.unige.ch}
\begin{document}
\maketitle
\begin{abstract}
Besides its construction as a quotient of an abelian surface, a Kummer surface can be obtained as the quotient of a K3 surface by a $\mathbb Z/2\mathbb Z $-action. In this paper, we classify all such K3 surfaces. 
Our classification is expressed in terms of period lattices and extends Morrison's criterion of K3 surfaces with a Shioda-Inose structure. Moreover, we list all the K3 surfaces associated to a general Kummer surface and provide very geometrical examples of this phenomenon. 
\end{abstract}

\section{Introduction}
The definition of a Kummer surface unravels the natural link between K3 surfaces and abelian surfaces. Indeed, let $A$ denote an abelian surface and let $A \stackrel{i} \to A$ be its involution automorphism, i.e. $i(a)=-a$. The quotient, $A/\{ 1, i \}$, is a singular surface with sixteen ordinary double points. Its minimum resolution, $\kum(A)$, is a K3 surface called the Kummer surface associated to $A$. Closely related to Kummer surfaces, are the K3 surfaces with a Shioda-Inose structure (\cite{SI}, \cite{Mor1}). A K3 surface $X$ is said to have a Shioda-Inose structure if it admits a rational map of degree two into a Kummer surface, $X \stackrel{\tau} \dashrightarrow \kum (A)$, which
induces  
a Hodge isometry of the period lattices,
$$\mr T_X (2) \stackrel{\tau_*} \to \mr T_{\kum (A)}.$$
Morrison proved that $X$ admits a Shioda-Inose structure if and only if there exists a Hodge isometry, $T_X \simeq T_A$, between the period lattices of $X$ and an abelian surface $A$ (\cite{Mor1}).

However, in general, a rational map of degree two between a K3 surface and a Kummer surface, does not need to induce a Hodge isometry of period lattices.

The aim of this paper is to treat this more general situation, i.e. it classifies all the K3 surfaces that admit a rational map of degree two into a given Kummer surface. As a result, we prove that if $X$ is a K3 surface such that its period lattice admits an embedding into the period lattice of an abelian surface $A$, $T_X \hookrightarrow T_A $, with the property that $$\quad T_A/T_X\simeq (\mathbb Z / 2)^ \alpha  \textrm{ with }0 \leq \alpha \leq 4,$$
then there exists a rational map of degree two, $$X \stackrel{\tau} \dashrightarrow \kum(A).$$ Note that the case where $\alpha =0$ recovers Morrison's criterion for K3 surfaces with a Shioda-Inose structure. Moreover, we give a complete list of all the K3 surfaces covering a general Kummer surface.

In the second part of this paper, we introduce the notion of an \textit{even eight} on a K3 surface. We use this concept to construct explicit examples of the surfaces appearing in our list. 
Interestingly in these cases, the rational map $X \stackrel{\tau} \dashrightarrow \kum(A)$ decomposes as
$$\xymatrix{X\ar@{-->}[d]_{\tau}
\ar[r]^{\zeta} & \Omega \ar[d]_{\varphi}\\
              \kum(A) \ar[r]^{\rho}& \mathbb P^2}$$
where $\rho$ is the double cover of $\mathbb P^2$ branched along a reducible sextic $\Sigma$ decomposable as a (possibly degenerate) quartic $\Gamma$ and a conic $\Delta$ and where $\varphi$ and $\zeta$ are double covers branched along $\Delta$ and $\varphi^*(\Gamma)$ respectively. 
\begin{acknowledgements}
  \begin{flushleft}
  This paper is part of the author's thesis work. The author would like to thank her adviser, Igor Dolgachev for all his guidance, suggestions and encouragement. The author is also very grateful to Bert van Geemen for his precious help in establishing the list of K3 surfaces covering a general Kummer surface.
  \end{flushleft}
\end{acknowledgements}

\section{Lattices and involutions on K3 surfaces}
The theory of K3 surfaces is intimately connected to the study of integral quadratic forms or \textit{lattices}. In this section, we recall the basic properties of lattices and discuss their fruitful applications to involutions on K3 surfaces.
A (non-degenerate) lattice $S$ is a free $\mathbb Z$-module of finite rank together with a (non-degenerate) bilinear form, $b: S \times S \to \mathbb Z$. The lattice $S$ is even if the associated quadratic form takes only even values. An embedding of lattices $M \hookrightarrow S$ is an homomorphism of $\mathbb Z$-modules preserving the bilinear forms. We say that the embedding is primitive if the quotient $S/M$ is torsion free.  The signature of $S$, $(t_+, t_-)$ is the signature of the bilinear form $b$, over the real vector space $S \otimes \mathbb R$. The lattice $S(m)$ for $m\in \mathbb Z $ denotes the same free $\mathbb Z$-module as $S$, with the bilinear form $b_{mS}(x,y)=mb_S(x,y).$
If $S$ is an non-degenerate even lattice, then the $\mathbb Q$-valued quadratic form on $S^*=\mathrm{Hom}(S,\mathbb Z)$ induces a quadratic form on the finite abelian group $D_S=S^*/S,$ $$q_S: D_S \to
\mathbb Q/ 2 \mathbb Z,$$ which is well defined mod $2\mathbb Z$ and a bilinear form $$b_S: D_S \times D_S
\to \mathbb Q/ \mathbb Z .$$ The pair $(D_S,q_S)$ is called the \textit{discriminant
form} of $S$, and $l(D_S)$ is the minimum number of generators of $D_S$. A non-degenerate lattice $S$ is unimodular if $D_S \simeq \{ 0 \}$. An embedding $S \hookrightarrow S'$ of even lattices with a finite
cokernel is called an \textit{overlattice}. Let $S \hookrightarrow S'$ be an overlattice of even non-degenerate
lattices, then $$|D_S|=[S':S]^2|D_{S'}|.$$
Consider the chain of embeddings $S \hookrightarrow S'
\hookrightarrow S'^* \hookrightarrow S^*$ and let $H_{S'}=S'/S$ be the subgroup of $D_S$. Then the correspondence, $S' \leftrightarrow H_{S'}$, determines a bijection
between even overlattices of $S$ and isotropic subgroups of $D_S$.


A $p$-adic lattice and its discriminant form are defined in the same way as for an integral lattice, only replacing $\mathbb Z$ by $\mathbb Z_p$. Let $S$ be a lattice over $\mathbb Z$, for every prime number $p$, $q_p$ denotes the restriction of the discriminant form $q_S$ to $(D_S)_p$, the $p$-component of $D_S$. There exists an unique $p$-adic lattice $K(q_p)$ whose discriminant form is isomorphic to $((D_S)_p, q_p)$ and whose rank is equal to $l((D_S)_p)$. The $p$-adic number $|K(q_p)|$ is then the determinant of a symmetric matrix associated to the lattice $K(q_p)$; it is therefore well defined mod $(\mathbb Z^*_p)^2$.

\begin{theorem}\label{primitive embedding} \cite{N1}
There exists a primitive embedding of the even lattice $S$, with invariants $(t_+, t_-, q_S)$ into the unimodular lattice  $L$ with signature $(l_+, l_-)$ if and only if the following conditions are simultaneously satisfied:
\begin{enumerate}
\item $l_+-l_- \equiv 0$ (mod 8),
\item $l_--t_- \geq 0$, $l_+-t_+ \geq 0$, $l_++l_- -t_+-t_-\geq l(D_S)$,
\item $(-1)^{l_+ -t_+}|D_S| \equiv \pm |K(q_p)| (\textrm{mod } \mathbb Z^*_p)^2$\\
for all odd primes p for which $ l_++l_- -t_+-t_-= l((D_S)_p),$
\item $|D_S| \equiv \pm |K(q_2)| (\textrm{mod }\mathbb Z^*_2)^2$\\
if $ l_++l_- -t_+-t_-= l((D_S)_2)$ and $q_2 \ncong \frac{x^2}{\theta \cdot 2} \oplus q_2'$ for some $\theta \in \mathbb Z^*_2, \textrm{ and some }q_2'.$
\end{enumerate}
\end{theorem}

A K3 surface $X$ is a smooth projective complexe surface satisfying

$$K_X=0 \quad \textrm{and} \quad \mr H^1(X, \mathcal O_X)=0.$$

\noindent Let $X$ be a K3 surface. The second cohomology group, $\mr H^2(X, \mathbb Z)$, equipped with the cup product, is an even unimodular lattice of signature $(3,19)$ and is isomorphic to $U^3 \oplus E_8^2(-1)$, where $U$ is the hyperbolic plane and $E_8(-1)$ is the unique even negative definite lattice of rank eight. The period lattice of a surface $X$, denoted by $T_X$, is defined by $$T_X = S_X ^\perp \subset \mr H^2(\mr X, \mathbb Z)$$ where $S_X$ is the N\'eron-S\'ev\'eri group of $X$.
Recall that the lattice $\mr H^2(X, \mathbb Z)$ admits a Hodge decomposition
of weight two $$\mr H^2(X, \mathbb C) \simeq \mr H^{2,0} \oplus \mr H^{1,1} \oplus \mr H^{0,2}.$$
Similarly, the period lattice $T_X$ inherits a Hodge
decomposition of weight two denoted by $$T_X \otimes \mathbb C
\simeq T^{2,0} \oplus T^{1,1} \oplus T^{0,2}.$$

\noindent An isomorphism between two lattices that preserves their bilinear forms and their Hodge decomposition is called a Hodge isometry. K3 surfaces have the unusual property that they satisfy a "`Torelli theorem"', namely

\begin{theorem}\label{torelli}\cite{torelli}
Let $X$ and $X'$ be two K3 surfaces and let $\phi: \mr H^2(X', \mathbb Z) \to  \mr H^2(X, \mathbb Z)$ be a Hodge isometry. If $\phi$ preserves effective classes, then there exists a unique isomorphism $u: X \to X'$ such that $u^*=\phi.$
\end{theorem}

\noindent Mukai improved this result for K3 surfaces of large picard number.

\begin{theorem}\label{M1} \cite{M1}
Let $X$ and $X'$ be K3 surfaces and let $\phi:T_X' \to T_X$ be a Hodge isometry. If $\rho(X) \geq 12$, there exists an isomorphism $u: X \to X'$ such that $\phi$ is induced by $u^*.$
\end{theorem}


\begin{definition}
Let $X$ be a K3 surface, an involution $X \stackrel{\iota} \to X$ is called \textit{a Nikulin involution} if it preserves any holomorphic two-form on $X$, i.e. $\iota^*(\omega)=\omega$ for all $\omega \in \mr H^{2,0}(X).$
\end{definition}

\begin{theorem}\cite{N5}
Let $X \stackrel{\iota} \to X$ be a Nikulin involution on $X$. Then

1) $\iota$ has eight isolated fixed points and the surface $X / \{ 1, \iota \}$ has eight ordinary double points,

2) the minimal resolution $Y$ of $X/ \{1, \iota \}$ is a K3 surface.
\end{theorem}

Let $\mr Z \stackrel{\epsilon} \to X$ be the blow-up of the eight fixed points of $\iota$. The involution $\iota$ lifts to an involution $\hat{\iota}$ on $Z$ whose fixed points consist of the eight exceptional curves. Thus the quotient surface, $Z \stackrel{p} \to Z / \{ 1, \hat{\iota}\} =Y$, contains eight disjoint smooth rational curves, $C_1, \dots, C_8$, which make up the branch locus of the map $p$. Consequently they satisfy $$C_1 + \dots + C_8 \in 2S_Y.$$ The commutative resolution diagram

$$\xymatrix{\mr Z \ar[r]^{\epsilon}\ar[d]_{p} & X\ar[d]_{/ \iota}  \\
 Y \ar[r]& X/ \{ 1, \iota \}} $$
induces a rational map of degree two between K3 surfaces, $X \stackrel{\tau} \dashrightarrow Y$.

\begin{definition}\label{Nikulin}
The Nikulin lattice is an even lattice $N$ of rank eight generated by $\{ c_i\}_{i=1}^8$ and $d=\frac{1}{2}\sum_{i=1}^8c_i$, with the bilinear form $c_i \cdot c_j = - 2\delta_{ij}.$
\end{definition}




Conversely, any rational map of degree two between K3 surfaces, $X \stackrel{\tau} \dashrightarrow Y$, is obtained from a similar construction and therefore it implies the existence of a Nikulin involution on $X$.




Let $X \stackrel{\iota} \to X$ be a Nikulin involution. The two invariant lattices $$T^{\iota}= \{x\in \mr H^2(X, \mathbb Z)|
\iota^*(x)=x\} \quad \textrm{and} \quad S_{\iota}= \{x\in \mr H^2(X, \mathbb Z)| \iota^*(x)=-x\}$$ define the primitive embeddings $$T_X \hookrightarrow T^{\iota}
\hookrightarrow \mr H^2(X, \mathbb Z), \quad S_{\iota}
\hookrightarrow S_X \hookrightarrow \mr H^2(X, \mathbb Z)$$
and the isomorphisms $$T^{\iota} \simeq U^3 \oplus E_8(-2)\quad
\textrm{and} \quad S_{\iota}\simeq E_8(-2).$$
\noindent Let $X \stackrel{\tau} \dashrightarrow Y$ be the rational map associated to $\iota$. The maps induced on the period lattices $$T_Y \stackrel{\tau^*}\to T_X \quad \textrm{and} \quad  T_X \stackrel{\tau_*}\to T_Y$$ preserve the Hodge decompositions and satisfy
$$\tau_*\tau^*(x)=2x, \quad \tau^*\tau_*(x)=2x, \quad \tau^*(x) \cdot \tau^*(y)=2xy, \quad \tau_*(x) \cdot \tau_*(y)=2xy.$$
\noindent The image of $T_Y$ by $\tau^*$ is then a sublattice of $T_X$ isomorphic to $T_Y(2)$ with two-elementary quotient, i.e. $$T_X / T_Y(2) \simeq (\mathbb
Z/ 2\mathbb Z)^\beta  \textrm{ for some } \beta \geq 0.$$
\begin{lemma}\label{Inv2}\cite{N3} Let $X$ be a $K3$
surface and let $T_X \hookrightarrow T \simeq U^3\oplus E_8(-2)$ be a primitive
embedding of lattices. Then there exists a Nikulin $\iota$ of $X$  such
that for the corresponding rational map of degree two $X \stackrel{\tau} \dashrightarrow Y$ $$\tau^*T_Y\simeq T_Y(2)\simeq2 (T^* \cap
(T_X \otimes \mathbb Q)) \subset T_X.$$
\end{lemma}


Finally, the following argument is due to Nikulin (\cite{N3}) and will be crucial for the proof of the theorem classifying the K3 surfaces covering a given Kummer surface. Let $X$ be a K3 surface with a Nikulin involution $\iota$, denote by $T$ the period lattice $T_X$ and by $S$, the invariant lattice $S_{\iota}$. Let $M$ be the primitive sublattice of $\mr H^2(X, \mathbb Z)$ generated by $T \oplus S$. Then the subgroup $\Gamma = M/(T \oplus S) \subset D_T \oplus D_S $ is isotropic and satisfies $$\Gamma \cap (D_T \oplus \{ 0 \})=0 \quad \textrm{and} \quad \Gamma \cap (\{ 0 \} \oplus D_S) =0.$$ Let $\pi_T$ and $\pi_S$ be the projections in $D_T$ and $D_S$ respectively. Let $\mathfrak H = \pi_T(\Gamma) \subset D_T$. If $\alpha= l(\mathfrak H) \leq 4$, then the discriminant group of $M$ is of the form $D_M \simeq D_T \oplus u(2)^{4-\alpha}$, where $u(2)$ denotes the discriminant form of the lattice $U(2).$ Note also that $\mathfrak{H}$ is not necessarily isotropic. The map $\xi= \pi_S \circ \pi^{-1}_T$ defines an inclusion of quadratic forms (i.e. an injection of groups preserving the bilinear form) $$\xi: \mathfrak H \to (D_S, -q_{S}).$$ Moreover, Nikulin proved the equality
$$((S^{\perp})^* \cap (T\otimes \mathbb Q))/T = \mathfrak H \subset D_T.$$


\section{Classification Theorem}

In this section, we give necessary and sufficient conditions for a K3 surface to admit a rational map of degree two into a given Kummer surface. As a corollary, we show that there exist finitely many such K3 surfaces and give an exhaustive list of all the K3 surfaces covering a general Kummer surface. 

\begin{theorem}\label{main thm}
Let $\mr{Kum}(A)$ be a Kummer surface. A K3 surface $X$ admits a rational map of degree two, $X
\stackrel{\tau} \dashrightarrow \kum(A)$, if and only if there exists an embedding
of period lattices $\phi: T_X \hookrightarrow T_A$ such that
\begin{enumerate}
\item $\phi_{\otimes \mathbb C}$ preserves the Hodge decomposition,
\item $T_A/T_X\simeq (\mathbb Z / 2 \mathbb Z)^ \alpha$ with $0 \leq \alpha \leq 4,$
\item $|D_{T_X}| \equiv \pm |K(q_{T_X})_2|$ mod ($\mathbb Z_2^*)^2$,\\
if $\alpha = \frac{\mathrm{rank}(T_X) + l((D_{T_X})_2)}{2}-3 \textrm{ and
} (q_{T_X})_2 \oplus u^{4-\alpha}(2) \ncong \frac{x^2}{\theta \cdot 2} \oplus q_2'$ for some $\theta \in \mathbb Z^*_2$ and some $q_2'$.
\end{enumerate}
\end{theorem}

\begin{proof}
Let $Y$ denote the Kummer surface $\kum(A)$ and suppose that $X$ admits a rational map of degree two into $Y$,
$X \stackrel{\tau} \dashrightarrow Y$. The covering involution $\iota: X \to X$ is then a Nikulin involution. Consider the induced Hodge isometry of lemma \ref{Inv2},
$$T_Y(2) \stackrel{\simeq} \to 2 (T_X \otimes \mathbb Q \cap (T^{\iota})^*) \subset T_X.$$
Denote by $\tilde{\mathfrak H}$ the rational overlattice of $T_X$ defined by $T_X \otimes \mathbb Q \cap (T^{\iota})^*$. Since $2\tilde{\mathfrak H} \subset T_X$, the quotient $\tilde{\mathfrak H}/ \ T_X$ is two-elementary. Since the period lattice of $\mr{Kum}(A)$ is Hodge isometric to $T_A(2)$ (\cite{N2}), $T_A(4)$ is Hodge isometric to $T_Y(2)$. Thus there is a Hodge isometry, $T_A(4) \stackrel{\simeq} \to 2\tilde {\mathfrak H} \subset T_X.$ The obvious identification of $T_A(4)$ with $2T_A$ shows that $T_A \simeq \tilde{\mathfrak H}$ and that the inclusion preserving the Hodge decomposition $T_X \hookrightarrow \tilde{\mathfrak H}\simeq T_A$ satisfies $$T_A/T_X \simeq \tilde{\mathfrak H}/ T_X \simeq  (\mathbb Z /2 \mathbb Z)^\alpha \textrm{ for some }\alpha \geq 0.$$

\noindent Moreover, the lattice $\tilde{\mathfrak H}$ must be an integral overlattice of $T_X$, and conseqently the quotient $\tilde {\mathfrak H}/ T_X$, is an isotropic subgroup of $(D_{T_X}, q_{T_X})$. Note that the group $\tilde {\mathfrak H}/ T_X$  corresponds to the group $\mathfrak H$ introduced in the last part of the previous section. Therefore $\mathfrak H$ admits an inclusion of quadratic forms into $(D_{S_{\iota}}, - q_{|S_{\iota}})$, where $l(D_{S_{\iota}})=8$. Because $\mathfrak H$ is isotropic, we obtain the upper bound $\alpha \leq 4$. Finally, let $M$ be the primitive closure of $T_X \oplus S_{\iota}$ in
$\mr H^2({X, \mathbb Z})$, then $D_{M} \simeq D_{T_X} \oplus u(2)
^{4- \alpha}$. The lattice $M$ admits a primitive embedding into the
unimodular lattice in $\mr H^2{(X, \mathbb Z)}$, so it satisfies the
hypothesis of the theorem \ref{primitive embedding}. In particular, it satisfies
 $$|D_{T_X}| \equiv \pm |K(q_{T_X})_2| \textrm{ mod }(\mathbb Z_2^*)^2$$
if $\alpha = \frac{\mathrm{rank}(T_X) + l(D_{T_X})_2}{2}-3 \textrm{ and
} (q_{T_X})_2 \oplus u^{4-\alpha}(2) \ncong
\frac{x^2}{\theta \cdot 2} \oplus q_2' \textrm{ for some } \theta \in \mathbb Z^*_2.$

Conversely, assume that $X$ satisfies conditions $(1)$, $(2)$ and $(3)$.

\noindent \textit{Claim:} There exists a primitive embedding $T_X
\hookrightarrow U^3 \oplus E_{8}(-2)$ such that $$(T_X \otimes
\mathbb Q)  \cap (U^3 \oplus E_8(-2)^*)\simeq T_A.$$ Assuming
the claim, by lemma \ref{Inv2}, there exists a Nikulin involution on $X$ such that the corresponding rational map of degree two $X \stackrel{\tau}
\dashrightarrow \tilde{Y}$ satisfies $$T_{\tilde{Y}}(2)
\simeq 2 (U^3 \oplus E_8(-2)^* \cap T_X \otimes \mathbb Q) \simeq 2T_A
\subset T_X.$$ The Hodge isometry $T_{\tilde{Y}}(2)
\simeq 2T_A $ implies that $T_{\tilde Y} \simeq T_A(2)$ and hence that $T_Y \simeq T_{\tilde Y}$. According to theorem \ref{M1}, $\tilde{Y} \simeq Y$.

\textit{Proof of the claim:} Let $\mathfrak H \simeq
(\mathbb Z/2 \mathbb Z)^{\alpha}$ be the isotropic subgroup of $D_{T_X}$
corresponding to the overlattice $T_X \subset T_A$. Since $\alpha \leq 4$, there
exists an embedding of quadratic forms $$\mathfrak H \stackrel{\xi}
\hookrightarrow (D_{S}, -q_{S}),$$ where $S \simeq E_8(-2)$. Consider
the subgroup $$\Gamma_{\xi} = \{ h+ \xi (h) | h \in \mathfrak H \}
\subset D_{T_X} \oplus D_{S}.$$ Clearly  $\Gamma_{\xi}$ is an
isotropic subgroup of $D_{T_X} \oplus D_{S}$. Denote by $M$ the
overlattice it defines and observe that $$\mr{rank}(T_M)=\mr{rank}(T_X) + 8 \textrm{ and } D_{M} \simeq
q_{T_X} \oplus u(2)^{4-\alpha}.$$ In order for $M$ to admit a primitive embedding into the unimodular lattice $\mr H^2(X, \mathbb Z)$, it must satisfy the hypothesis of the theorem \ref{primitive embedding}. In our setting it involves checking that $$ l((D_M)_p) + \mathrm{rank}(T_X) \leq 14 \quad \textrm{for all } p \textrm{ prime},$$since for any finite abelian group, $l(D)= \mathrm{max}_{p}l(D_p)$. If $p$ is odd, then $(D_M)_p \simeq (D_{T_X})_p$ and the strict inequality $l((D_{T_X})_p) + \mathrm{rank}(T_X) < 14$ always holds, as $l((D_{T_X})_p)\leq \mr{rank}(T_X)$ and $2 \leq \mr{rank}(T_X)\leq 5$. If $p$ is $2$, then $l((D_M)_2) = l((D_{T_X})_2)+ 8 -2\alpha$. The facts that  $l((D_{T_X})_2)-2\alpha\leq l((D_A)_2)$ (see the technical lemma \ref{technical lemma} below), and that the lattice $T_A$ embeds primitively into $U^3$ imply the inequalities $$l((D_{T_X})_2)+ 8 -2\alpha + \mathrm{rank}(T_X)\leq l((D_{T_A})_2)+8 +\mathrm{rank}(T_A) \leq l(D_{T_A}) + 8+ \mathrm{rank}(T_A)\leq 6+8 =14.$$
Thus the inequality $l((D_M)_p) + \mathrm{rank}(T_X) \leq 14 $ holds for all $p$.   The assumption
$$|D_{T_X}| \equiv \pm |K(q_{T_X})_2| \textrm{ mod }(\mathbb Z_2^*)^2$$
if $\alpha = \frac{\mathrm{rank}(T_X) + l(D_{T_X})_2}{2}-3 \textrm{ and
} (q_{T_X})_2 \oplus u^{4-\alpha}(2) \ncong
\frac{x^2}{\theta \cdot 2} \oplus q_2' \textrm{ for some } \theta \in \mathbb Z^*_2$, is therefore enough to guarantee the
existence of the primitive embedding of $M$ into $\mr H^2(X, \mathbb Z)$. Consider the embedding $$T_X \oplus S \subset M \subset
\mr H^2(X, \mathbb Z).$$ By construction, its restriction to
$T_X$ and to $S$ are primitive. It follows that $S$ is primitively
embedded into $S_X$ or equivalently $T_X$ is primitively embedded into
$S^{\perp} \simeq U^3 \oplus E_8(-2).$ Finally, from the last equality of section 2, we find that $$\mathfrak H \simeq (T_X \otimes \mathbb Q)  \cap (U^3 \oplus
E_8(-2)^*)/ T_X$$ and therefore $$(T_X \otimes \mathbb Q)  \cap
(U^3 \oplus E_8(-2)^*)\simeq T_A.$$
\end{proof}

\begin{lemma}\label{technical lemma}
Let $D$ be a finite abelian $2$-group with a bilinear form, $b(x,y) \in \mathbb Q/ \mathbb Z$ and a quadratic form $q(x)=b(x,x) \in \mathbb Q/ 2\mathbb Z$. Let $\mathfrak H \subset D$ be an isotropic subroup, then
$$l(D)- 2l(\mathfrak H) \leq l(\mathfrak H^{\perp}/ \mathfrak H).$$
\end{lemma}
\begin{proof}
It is enough to show that $l(D/\mathfrak H) \leq l(\mathfrak H^{\perp})$.
Indeed it would imply that $$l(D) \leq l(D/\mathfrak H) + l(\mathfrak H) \leq l(\mathfrak H^{\perp}) + l(\mathfrak H) \leq  l(\mathfrak H^{\perp}/ \mathfrak H) + l(\mathfrak H) + l(\mathfrak H)$$
where the first and last inequalities are obvious.

Assume first that $\mathfrak H \simeq \mathbb Z /2$ with generator $h$.
Let $[x] \in D/\mathfrak H$ be a non-zero element of order $2^r$, then
$$q(h)=b(h,h)=b(2^rx,h)=2b(2^{r-1}x,h)=0  \Rightarrow b(2^{r-1}x,h)=0 \textrm{ therefore } 2^{r-1}x \in \mathfrak H^{\perp}.$$
Under this correspondence, a minimal set of generators of $D/\mathfrak H$ maps to a set of independent elements of $\mathfrak H^{\perp}$, i.e. $l(D/\mathfrak H) \leq l(\mathfrak H^{\perp})$.

If $|\mathfrak H|>2$, then there exists a subgroup $\mathfrak H' \subset \mathfrak H$, such that $\mathfrak H/\mathfrak H' \simeq \mathbb Z/2$. Clearly $\mathfrak H/\mathfrak H'$ is an isotropic subgroup of $D_S/\mathfrak H'$.  Since a minimal set of generators of $(\mathfrak H/\mathfrak H')^{\perp}$ define a set of independent elements of $\mathfrak H^{\perp}$, we get the inequality $l((\mathfrak H/\mathfrak H')^{\perp}) \leq l(\mathfrak H^{\perp})$.
By induction, we obtain $$l(D/\mathfrak H) =l((D/\mathfrak H')/\mathfrak H/\mathfrak H') \leq l((\mathfrak H/\mathfrak H')^{\perp}) \leq l(\mathfrak H^{\perp}).$$
\end{proof}

Let $A$ be an abelian surface. If $X$ is a K3 surface such that its period lattice $T_X$ admits an embedding into $T_A$ with $T_A/T_X \simeq (\mathbb Z/2\mathbb Z)^{\alpha}$, then there is an embedding of $2T_A$ into $T_X$. The finite group $T_X/2T_A$ is a subspace of the $\mathbb F_2$-symplectic space $(T_A /2T_A, q_{T_A/2T_A})$. Conversely, any subspace $U$ of $T/2T_A$ defines a sublattice of $T_A$, $T_U$, containing $2T_A$ and hence satisfying $T_A/T_U \simeq (\mathbb Z/2\mathbb Z)^{\alpha}$.  Moreover, the lattice $T_U$ inherits a Hodge decomposition from the one on $T_A$ and if $T_U$ satsifies the condition $(3)$ of the theorem \ref{main thm}, it determines an unique K3 surface. Note that any automorphism of $T_A$ preserving its Hodge decomposition must be $\pm$ identity,\cite{torelli}, and therefore if $U$ and $U'$ are two subspaces of $T_A/2T_A$ such that there exists an orthogonal transformation $\phi \in O(T_A/2T_A)$ mapping $U$ to $U'$, then (if it exists) the lifting of $\phi$ to $T_A$ cannot define a Hodge isometry between $T_U$ and $T_{U'}$. Clearly, we have the following corollary.
\begin{corollary}
Let $\kum(A)$ be a Kummer surface. There exist finitely many
isomorphism classes of K3 surfaces admitting a rational map of
degree two into $\kum(A)$.
\end{corollary}

We say that a Kummer surface $Y$ is general if it is associated to a general abelian surface $A$, that is to say $A$ is the jacobian of a curve of genus two and its period lattice $T_A$ is isomorphic to $U \oplus U \oplus \langle -2 \rangle$. In order to determine all the sublattices of $T_A$ satisfying the conditions of \ref{main thm}, we classify all the subspaces of $(T_A/ 2T_A \simeq \mathbb F_2^5, q(x)=x_1x_2+ x_3x_4 +x^2_5$).
%
We obtain in this way a short and exhaustive list of all K3 surfaces satisfying the conditions of theorem \ref{main thm}. In total, we find 373 K3 surfaces that are listed below. The number in the first column indicates the number of non-isomorphic K3 surfaces within a given type of period lattice. In other words, it represents the size of the orbit of a given subspace of $T_A/2T_A$ for the action of the orthogonal group $O(\mathbb F_2^5, q)$.


\vspace{6pt}
\begin{center}
\begin{minipage}{\textwidth}
\begin{center}
{\fontsize{10}{16}\selectfont\begin{tabular}{cccc}

\hline
$\sharp$ K3's& $T_X$ & $D_{T_X}$ &
$q_{T_X}$ \\

\hline $1$ & $U \oplus U \oplus \langle-2\rangle$ &
$\mathbb Z / 2 \mathbb Z$ & $\frac{1}{2}x^2$\\

\hline $15$ & $U(2) \oplus U \oplus \langle-2\rangle $ & $\mathbb Z / 2 \mathbb Z \oplus \mathbb Z / 2 \mathbb Z \oplus\mathbb Z / 2 \mathbb Z $ & $x_1x_2 - \frac{1}{2}x_3^2$ \\

\hline $10$ & $U \oplus U \oplus \langle-8\rangle$ & $\mathbb Z / 8 \mathbb Z$ & $-\frac{1}{8}x^2$ \\

\hline $6$ & $U \oplus N$ & $\mathbb Z / 8 \mathbb Z$ & $\frac{3}{8}x^2$ \\

\hline $20$ & $U(4) \oplus U \oplus \langle-2\rangle$ & $\mathbb Z / 4 \mathbb Z \oplus \mathbb Z / 4 \mathbb Z \oplus \mathbb Z / 2 \mathbb Z$ & $\frac{1}{2}x_1x_2- \frac{1}{2}x_3^2$ \\

\hline $15$ & $U(2) \oplus U(2) \oplus \langle-2\rangle$ & $\mathbb Z / 2 \mathbb Z ^{\oplus 5}$ & $x_1x_2 + x_3x_4 - \frac{1}{2}x_5^2$ \\

\hline $45$ & $U(2) \oplus U \oplus \langle-8\rangle$ & $\mathbb Z / 2 \mathbb Z \oplus \mathbb Z / 2 \mathbb Z \oplus \mathbb Z / 8 \mathbb Z $ & $x_1x_2  - \frac{1}{8}x_3^2$ \\

    \hline $60$ & $U \oplus \langle -2 \rangle \oplus  \langle 2 \rangle \oplus  \langle -8 \rangle$ & $\mathbb Z / 2 \mathbb Z \oplus \mathbb Z / 2 \mathbb Z \oplus \mathbb Z / 8 \mathbb Z $ & $-\frac{1}{2}x_1^2+\frac{1}{2}x_2^2 -\frac{1}{8}x_3^2 $ \\
\hline $15$ & $U(2) \oplus N$ & $\mathbb Z / 2 \mathbb Z \oplus \mathbb Z / 2 \mathbb Z \oplus \mathbb Z / 8 \mathbb Z $ & $x_1x_2 +\frac{3}{8}x_3^2$ \\
\hline $15$ & $U(4) \oplus U(2) \oplus \langle-2\rangle$ & $  \mathbb Z / 4 \mathbb Z ^{\oplus 2} \oplus \mathbb Z / 2 \mathbb Z ^{\oplus 3}$ & $\frac{1}{2}x_1x_2 + x_3x_4 - \frac{1}{2}x_5^2$ \\

\hline $60$ & $U(4) \oplus U \oplus \langle-8\rangle$ & $\mathbb Z / 4 \mathbb Z \oplus \mathbb Z / 4 \mathbb Z \oplus \mathbb Z / 8 \mathbb Z $ & $\frac{1}{2}x_1x_2  - \frac{1}{8}x_3^2$ \\

\hline $20$ & $U(4) \oplus N$ & $\mathbb Z / 4 \mathbb Z \oplus \mathbb Z / 4 \mathbb Z \oplus \mathbb Z / 8 \mathbb Z $ & $\frac{1}{2}x_1x_2 -\frac{1}{2}x_1^2 + x_1x_3 + \frac{3}{8}x_3^2$ \\

\hline $15$ & $U(2) \oplus U(2) \oplus \langle-8\rangle$ & $ \mathbb Z / 2 \mathbb Z ^{\oplus 4} \oplus \mathbb Z / 8 \mathbb Z  $ & $x_1x_2 + x_3x_4 - \frac{1}{8}x_5^2$ \\

\hline $45$ & $\langle 2 \rangle ^2 \oplus \langle -2 \rangle ^2 \oplus \langle -8 \rangle$ & $ \mathbb Z / 2 \mathbb Z ^{\oplus 4} \oplus \mathbb Z / 8 \mathbb Z  $ & $\frac{1}{2}(x_1^2+ x_2^2 -x_3^2 -x_4^2)- \frac{1}{8}x_5^2$ \\

\hline $1$ & $U(4) \oplus U(4) \oplus \langle-2\rangle$ & $\mathbb Z / 4 \mathbb Z ^{\oplus 4} \oplus \mathbb Z / 2 \mathbb Z $ & $\frac{1}{2}x_1x_2 + \frac{1}{2} x_3x_4 - \frac{1}{2}x_5^2$ \\

\hline $15$ & $U(4) \oplus U(2) \oplus \langle-8\rangle$ & $\mathbb Z / 2 \mathbb Z ^{\oplus 2} \oplus \mathbb Z / 4 \mathbb Z ^{\oplus 2} \oplus \mathbb Z / 8 \mathbb Z $ & $x_1x_2 + \frac{1}{2} x_3x_4 - \frac{1}{8}x_5^2$ \\

\hline $15$ & $U(4) \oplus \langle2\rangle \oplus \langle-2\rangle \oplus \langle-8\rangle$ & $\mathbb Z / 2 \mathbb Z ^{\oplus 2} \oplus \mathbb Z / 4 \mathbb Z ^{\oplus 2} \oplus \mathbb Z / 8 \mathbb Z $ & $\frac{1}{2}x_1^2 - \frac{1}{2}x_2^2 + \frac{1}{2}x_3x_4 - \frac{1}{8}x_5^2$ \\
\hline

\end{tabular}\label{list}}
\end{center}
\end{minipage}
\end{center}
\vspace{20pt}

\noindent The bilinear form of the lattices $N$ is given by the symmetric integral matrix
$$N=\left( \begin{smallmatrix} 2 & 1 & 2 \\   1 & -2 & 0 \\ 2& 0&  0\end{smallmatrix}\right).$$



\section{Even Eight}
The remaining of this paper will be devoted to the geometrical description of the K3 surfaces appearing in the above list that have two-elementary period lattices. The study of such K3 surfaces is made easier by the existence on them of another involution.

\begin{prop}\label{2-elem}\cite{N1}
Let $X$ be a K3 surface. The period lattice $T_X$ is two-elementary if and only if
$X$ admits an involution $\theta: X \to X$ such that
$\theta^*_{S_X}=\mr{id}_{S_X}$ and $\theta^*_{T_X}=-\mr{id}_{T_X}$.
\end{prop}

\noindent As a first step towards our geometrical description, we introduce the notion of an even eight. Secondly, in order to provide examples of different even eights on a Kummer surface, we give explicit generators for the N\'eron-S\'ev\'eri lattice of a general Kummer surface. We follow Naruki's notation in \cite{Naruki}.


\begin{definition}
Let $Y$ be a K3 surface, an even eight on $Y$ is a set of eight disjoint smooth rational curves, $C_1, \cdots, C_8$, such that $C_1 + \cdots + C_8 \in 2S_Y.$
\end{definition}

Recall, from section 2, that a rational map of degree two between K3 surfaces, $X \stackrel{\tau} \dashrightarrow Y$  factors as
$$\xymatrix{\mr Z \ar[rr]^{\epsilon}\ar[rd]_{p} && X\ar@{-->}[ld]^{\tau}  \\
 &Y&} $$
where $\epsilon$ is the blow-up of $X$ at eight points and $p$ is a regular map of degree two, branched along eight disjoint smooth rational curves, hence an even eight. Conversely, given an even eight, $C_1, \cdots, C_8$, on a K3 surface $Y$, there is exists a double cover $Z \stackrel{p} \to Y$ branched along $\sum_{i=1}^{8} C_{i}$. Let $p^*(C_i)=2E_i$, then each $E_i$ is an exceptional divisor of the first kind. We may therefore blow down $\sum_{i=1}^{8}
E_{i}$ to a smooth surface $X$, which is known to be a K3 surface. We recover, in this way, a rational map of degree two $ X \stackrel{\tau}\dashrightarrow Y$ between the K3 surfaces $X$ and $Y$ (see \cite{N2}).







Let $Y \simeq \kum(A)$ be a general Kummer surface and let $C$ be the genus two curve of which $A$ is the jacobian. The degree two map given by the linear system $|2C|$, $A \stackrel{|2C|} \to \mathbb P^3$, factors through the involution $i$, where $i(a)=-a$, and hence defines an embedding $A/\{ 1, i\} \hookrightarrow \mathbb P^3$. The image of this map, $Y_0$, is a quartic in $\mathbb P^3$ with sixteen nodes. Denote by $L_0$ the class of a hyperplane section of $Y_0$. Projecting from a node, $Y_0$ admits a regular map of degree two into the projective plane which induces the map $Y \stackrel{\phi} \to \mathbb P^2$, given by the linear system $|L-E_0|$, where $L$ is a the pullback of $L_0$ on $Y$ and $E_0$ is the exceptional curve resolving the center of projection. The branch locus of the map $\phi$ is a reducible plane sextic, which is the union of six lines, $l_1, \cdots, l_6$, all tangent to a conic $W$.

\begin{figure}[h]
$$
  \begin{xy}
   <0cm,0cm>;<1.5cm,0cm>:
    (2,-.3)*++!D\hbox{$l_1$},
    (1.1,0.25)*++!D\hbox{$l_2$},
    (1.1,1.3)*++!D\hbox{$l_3$},
    (2,1.8)*++!D\hbox{$l_4$},
    (2.8,1.3)*++!D\hbox{$l_5$},
    (2.9,0.2)*++!D\hbox{$l_6$},
    (1.5,1)*++!D\hbox{$W$},
    (.5,0.18);(3.5,0.18)**@{-},
    (0.5,2);(1.6,0)**@{-},
    (3.1,2);(2.6,0)**@{-},
    (0.7,.5);(2,2.5)**@{-},
    (0.5,1.82);(3.5,1.82)**@{-},
    (2,2.7);(3.45,0)**@{-},
    (2,1)*\xycircle(.8,.8){},
  \end{xy}
  $$
\caption{\label{plane sextic}}
\end{figure}

\noindent Let $p_{ij}=l_i \cap l_j \in \mathbb P^2$, where $1\leq i < j \leq 6$. Index the ten $(3,3)$-partitions of the set $\{ 1, 2, \dots ,6 \}$, by the pair $(i,j)$ with $2 \leq i < j \leq 6$. Each pair $(i,j)$ defines a plane conic $l_{ij}$ passing through the sixtuplet $p_{1i},
p_{1j}, p_{ij},p_{lm},p_{ln},p_{mn}$, where $\{ l,m,n \}$ is the complement of $\{ 1,i,j\}$ in $\{ 1, 2, \dots ,6 \}$ and where $l \leq m \leq n$. The map $\phi$ factors as $Y \stackrel{\tilde{\phi}} \longrightarrow \tilde{\mathbb P}^2\stackrel{\eta}\longrightarrow \mathbb P^2$ where $\eta$ is the blowup of
$\mathbb P^2$ at the $p_{ij}$'s and where $\tilde{\phi}$ is the double cover of $\tilde{\mathbb P}^2$ branched along the strict transform in $\tilde{\mathbb P}^2$ of the plane sextic of figure \ref{plane sextic}. Denote by $E_{ij} \subset Y$ the preimage of the exceptional curves of $\tilde{\mathbb P}^2$. The ramification of the map $\tilde{\phi}$ consists of the union of six disjoint smooth rational curves, $C_0+ C_{12} + C_{13}+ C_{14}+ C_{15}+ C_{16}$. The preimage of the ten plane conics $l_{ij}$ defines ten more smooth disjoint rational curves on $Y$, $C_{ij}, 2 \leq i < j \leq 6$. Finally, note that $\phi(E_0)=W$. The sixteen curves $E_0, E_{ij} \quad 2 \leq i < j \leq 6$ are called the nodes of $Y$ and the sixteen curves $C_0, C_{ij} \quad 2 \leq i < j \leq 6$ are called the tropes of $Y$. These two sets of smooth rational curves satisfy a beautiful configuration called the $16_6$-configuration, i.e. each node intersects exactly six tropes and vice versa.

\noindent It is now possible to fully describe the N\'eron-S\'ev\'eri lattice $S_{Y}$ of a general Kummer surface. Indeed, it is generated by the
classes of $E_{0}, E_{ij}$, $C_{0}, C_{ij}$
and $L$, with the relations:

\begin{enumerate}
\item $C_{0}= \frac{1}{2}(L- E_{0} -
\sum_{i=2}^{6} E_{1i}),$
\item $C_{1j}= \frac{1}{2}(L - E_{0} -
\sum_{i\neq j}E_{ij}),$ where $2\leq j \leq 6,$
\item $C_{jk}= \frac{1}{2}(L - E_{1j} -E_{1k}- E_{jk}- E_{lm}-E_{ln}-
E_{mn})$ where $2\leq i < j \leq 6,$ and $\{ l, m, n\}$ are as described above.

\end{enumerate}

\noindent The intersection pairing is given by:

\begin{enumerate}

\item  the $E_{0}, E_{ij}$ are mutually orthogonal,

\item $\langle L, L \rangle =4, \langle L, E_{0} \rangle =\langle L, E_{ij}\rangle =0,$

\item $\langle E_{0}, E_{0} \rangle= \langle E_{ij}, E_{ij} \rangle=-2$,

\item the $C_{0}, C_{ij}$ are mutually orthogonal,

\item  $\langle L, C_{0} \rangle=  \langle L, C_{ij} \rangle=2$.

\end{enumerate}


\noindent The action on $S_Y$ of the covering involution $\alpha$ of the map $\phi$ is given by:

\begin{tabular}{lrrlr}

$\alpha (C_{0})=C_{0}$ & &&

$\alpha (C_{1j})=C_{1j}$ & $2 \leq j \leq 6$\\

$\alpha (E_{ij})=E_{ij}$ & $1 \leq i < j \leq 6$ &&

$\alpha (L)=3 \mathrm L - 4 E_{0}$ & \\

$\alpha (E_{0})=2 L - 3E_{0}$& &&

$\alpha (C_{ij})= C_{ij} + L - 2E_{0}$ & $2 \leq i < j \leq 6$.\\

\end{tabular}

\begin{remark}
Conversely, the minimal resolution of the double cover of $\mathbb P^2$ branched along the sextic of figure \ref{plane sextic} is a Kummer surface (see \cite{H1} for a proof).\\
\end{remark}
\noindent Consider the divisors

\begin{tabular}{l}
$e_1:=(L-E_0)-(E_{12}+E_{46}),$ \\
$e_2:=2(L-E_0)-(E_{12}+E_{13}+E_{24}+E_{46}+E_{56}),$ \\
$e_3:=3(L-E_0)-2E_{12}-(E_{13}+E_{24}+E_{36}+E_{45}+E_{46}+E_{56}),$\\
$e_4:=4(L-E_0)-2(E_{12}+E_{13}+E_{46})-(E_{24}+E_{25}+E_{36}+E_{45}+E_{56}),$\\
$e_5:=5(L-E_0)-3E_{12}-2(E_{13}+E_{46}+E_{56})+(E_{24}+E_{25}+E_{34}+E_{36}+E_{45}),$\\
\end{tabular}

\begin{tabular}{ccc}
$e_6:=C_{23},$& $e_7:=\alpha(C_{23}),$ &$e_8:=E_{35}.$
\end{tabular}

\begin{theorem}\label{Naruki}\cite{Naruki}
The eight $(-2)$-classes, $e_1, \dots, e_8$ represent an even eight on $Y$.

\end{theorem}

\noindent The geometric meaning of the $(-2)$-classes can be explained in terms of the double plane model of the general Kummer surface, $Y \stackrel{\phi} \to \mathbb P^2$.
The class $e_1$ is represented by the proper inverse image of a line passing through $p_{12}$ and $p_{46}$. The class of $e_2$ is represented by the proper inverse image of a conic passing through $p_{12}, p_{13}, p_{24}, p_{46}, p_{56}$. The class of $e_3$ is represented by the proper inverse image of a cubic passing through $p_{13}, p_{24}, p_{36}, p_{45}, p_{46}, p_{56}$ and having a double point at $p_{12}$. The class of $e_4$ is represented by the proper inverse image of a quartic passing through $p_{24}, p_{25}, p_{36}, p_{45}, p_{56}$ and having double points at $p_{12}, p_{13}, p_{46}$. The class of $e_5$ is represented by the proper inverse image of a quintic passing through $p_{24}, p_{25}, p_{34}, p_{36}, p_{45}$ and having double points at $p_{13}, p_{46}, p_{56}$ and a triple point at $p_{12}$.


\section{K3 surface with a Shioda-Inose structure}
A particular case in the list of page 8 is given by the unique K3 surface $X_0$ that has its period lattice isomorphic to $U \oplus U \oplus \langle -2 \rangle$. In fact, Naruki showed that the map $X_0 \stackrel{\tau} \dashrightarrow Y$ is the rational map of degree two associated to the even eight of the theorem \ref{Naruki}. Moreover, according to Morrison's criterion, $X_0$ has a Shioda-Inose structure (\cite{Mor1}).

Prior to describing the geometry of $X_0$, we will state some basics about elliptic fibrations on K3 surfaces (see \cite{Mir1} for details).

\begin{definitions}
\begin{enumerate}

\item A smooth surface $X$ is \textit{elliptic} if it admits a
regular map into a curve $C$, $X \stackrel{f} \to C$ where the general fiber is
a smooth connected curve of genus one.

\item When nonempty, the set of sections of $f,$ $MW_f(X)$, form a
finitely generated group called the Mordell-Weil group, which can be identified with the set of
$k(C)$-rational points of the generic fiber $X_\eta$.

\item Let $F$ denote the class of a fiber. The set $\{ B \cdot F | B \in S_X \}$ is an ideal of $\mathbb Z$. The multisection index $l$ of
$f$ is the positive generator of this ideal. An $l-$section is a divisor $B$ on $X$ for which $B \cdot F =l$.
\end{enumerate}

\end{definitions}

\noindent We refer the reader to \cite{Mir1} for a table of Kodaira's classification of singular fibers on an elliptic fibration.

\begin{theorem}\label{MW}\cite{Shioda}
Let $X \stackrel{f} \to \mathbb P^1$ be an elliptic $K3$ surface with a section
$S_0$.  Let $A$ be the subgroup of $S_X$ generated by the section $S_0$,
the class of a fiber and the components of the reducible fibers
that do not meet $S_0$.
There is an short exact sequence of groups $$0 \to A
\stackrel{\alpha}\to S_X \stackrel{\beta}\to MW_f(X)\to 0 $$ where $\beta$ is the restriction to the generic
fiber. Let $\Theta_{i}$ $(i=1, \dots ,k)$ be all
the singular fibers of $f$. Then
\begin{enumerate}
\item $24=\chi_{top}(X)=\sum_{i}\chi_{top}(\Theta_{i}),$
\item $\mr{rank}(\mathcal
MW_{f}(X))=\mr{rank}(S_{X})-2-\sum_{i}(m(\Theta_{i})-1),$ where
$m(\Theta_{i})$ denotes the number of irreducible components of
$\Theta_{i}.$
\item When the rank of $\mathcal MW_{f}(X)$ is zero, then $$|D_{S_X}|=\frac{\prod_{i}s_{i}}{n^2}$$
where $s_{i}$ is the number of simple components of the singular
fiber $\Theta_{i}$ and $n$ is the order of the finite group
$\mathcal MW_{f}(X)$.
\end{enumerate}
\end{theorem}

\noindent Pjatecki{\u\i}-{\v{S}}apiro and {\v{S}}afarevi{\v{c}}
showed that on a K3 surface, if $D$ is a non-zero effective nef
divisor with $D^2=0$, then the system $\vert D \vert$ is of the
form $|mD_0|, m\geq 1$ and the linear system $|D_0|$ defines a
morphism $X \stackrel{|D_0|} \to \mathbb P^{1}$, with general
fiber a smooth connected curve of genus one (\cite{torelli}).


\begin{prop}\label{San1}
Let $Y$ be a general Kummer surface. Then there exists an elliptic fibration on $Y$ whose singular fibers are of the
type $6I_2 + I^*_5 + I_1$ and whose Mordell-Weil group is finite of order two.
\end{prop}
\begin{proof}
Consider the divisor class
$$D=5(L-E_0)-3E_{12}-2(E_{13}+E_{46}+E_{56})-(E_{24}+E_{25}+E_{36}+E_{45}).$$

\noindent Geometrically, $D$ can be represented by the proper inverse image, under the map $Y \stackrel{\phi} \to \mathbb P^2$, of a quintic passing trough $p_{24}$, $p_{25}$, $p_{36}$, $p_{45}$, having double points at $p_{13}$, $p_{46}$, $p_{56}$ and a triple point at $p_{12}$. Since  $D \sim e_5 + E_{34}$ and $e_5 \cdot E_{34}=2$, $D$ must be a reduced nef effective divisor with $D^2=0$. Let $Y \stackrel{f} \to \mathbb P^1$ be the elliptic fibration defined by the linear system $|D|$.
The divisors \begin{enumerate} \item $F_1=e_5+ E_{34}$
\item $ F_2=
e_4+ \underbrace{(L-E_0)-(E_{12}+E_{56})}_{a_1}$ \item $F_3=e_3+
\underbrace{2(L-E_0)-(E_{12}+E_{13}+E_{46}+E_{56}+E_{25})}_{a_2}$\item $F_4=e_2+
\underbrace{3(L-E_0)-2E_{12}-(E_{13}+E_{46}+E_{56}+E_{25}+E_{36}+E_{45})}_{a_3}$
\item $F_5=e_1+ \underbrace{4(L-E_0)-2(E_{12}+E_{13}+E_{56})-(E_{46}+E_{24}+E_{25}+E_{36}+E_{45})}_{a_4}$
\item $F_6=e_8+ \underbrace{5(L-E_0)-3E_{12}-2(E_{13}+E_{46}+E_{56})-(E_{24}+E_{25}+E_{36}+E_{45}+E_{35})}_{a_5}$
\end{enumerate}
represent six fibers of type $I_2$, as under the genericity assumption, the $(-2)$-cycles, $a_1, \dots a_{5}$ are represented by smooth rational curves. Their geometric meaning can be explained exactly in the same way as in the theorem \ref{Naruki}.
The divisor
$$F_7=e_6+e_7+2(E_{23}+C_{12}+E_{26}+C_{16}+E_{16}+C_0)+E_{14}+E_{15}$$
defines a fiber of type $I^*_5$. Note also that $C_{14} \cdot D= C_{15} \cdot D =1$, i.e. the curves $C_{14}$ and
$C_{15}$ are sections of $f$. The fifteen components of the singular fibers, $e_{1}, \dots, e_{8}, C_{12}, E_{23}, E_{26}, C_{16}, E_{16},$
$ C_0, E_{14},$ are independent over $\mathbb Q$. Consequently, by the theorem \ref{MW}, there are no other reducible fibers and the Mordell-Weil group is finite of order two.
\end{proof}

\noindent We note that the curves $e_1, \cdots, e_8$, defining the even eight of theorem \cite{Naruki} are components of the singular fibers of this elliptic fibration. In fact, we see that the surface $X_0$ admits also an elliptic fibration, $X_0 \stackrel{|\tau^*D|} \to \mathbb P^1$, with singular fibers of type $I_2+ I^*_{10}+6 I_2$ and with a Mordell-Weil group of order two. The reducible fibers and the sections of this latter fibration form the diagram of curves on $X_0$, shown in figure \ref{AX}. Nikulin proved that those curves are the only smooth rational curves lying on $X_0$ and that the group $\mr{Aut}(X_0)$ is generated by two commuting involutions $\theta$ and $\iota$. The involution $\theta$ satisfies $\theta^*_{S_{X_0}}=\mr{id}_{S_{X_0}}$, $\theta^*_{T_{X_0}}=-\mr{id}_{T_{X_0}}$, while the involution $\iota$ is the Nikulin involution that gives the Shioda-Inose structure on $X_0$ (\cite{N4}).
\begin{figure}[h]
$$ \begin{xy}
<0cm,0cm>;<1cm,0cm>:
    (0,0)*{\bullet},
    (0,0)*++!D\hbox{$S_{1}$},
    (0,0),(1,1)**@{-},
    (1,1)*{\bullet},
    (1,1);(2,1)**@{-},
    (1,1)*++!D\hbox{$R_{2}$},
    (2,1)*{\bullet},
    (2,1);(3,1)**@{-},
    (2,1)*++!D\hbox{$S_{2}$},
    (3,1)*{\bullet},
    (3,1);(4,1)**@{-},
    (3,1)*++!D\hbox{$R_{3}$},
    (4,1)*++!D\hbox{$S_{3}$},
    (4,1)*{\bullet},
    (4,1);(5,1)**@{-},
    (5,1)*++!D\hbox{$R_{4}$},
    (5,1)*{\bullet},
    (6,1);(5,1)**@{-},
    (6,2)*{\bullet},
    (6,2)*++!D\hbox{$R$},
    (5,1);(6,2)**@{-},
    (6,1)*++!D\hbox{$S_{4}$},
    (6,1)*{\bullet},
    (7,1)*{\bullet},
    (6,1);(7,1)**@{-};
    (7,-1);(6,-1)**@{-};
    (7,1)*++!D\hbox{$R_{5}$},
    (7,-1)*{\bullet},
    (7,-1)*++!D\hbox{$R'_{5}$},
    (8,1)*{\bullet},
    (8,1)*++!D\hbox{$S$},
    (7,1);(8,1)**@{-},
    (8,-1)*{\bullet},
    (8,-1)*++!U\hbox{$S'$},
    (7,-1);(8,-1)**@{-},
    (8,1);(8,-1)**@{=},
    (6,-1)*{\bullet},
    (6,-1);(5,-1)**@{-},
    (6,-1)*++!U\hbox{$S_{5}$},
    (5,-1)*{\bullet},
    (5,-1);(4,-1)**@{-},
    (6,-2)*{\bullet},
    (6,-2)*++!U\hbox{$T$},
    (5,-1);(6,-2)**@{-},
    (5,-1)*++!D\hbox{$R_{6}$},
    (4,-1)*{\bullet},
    (3,-1);(4,-1)**@{-},
    (4,-1)*++!U\hbox{$S_6$},
    (3,-1)*{\bullet},
    (2,-1);(3,-1)**@{-},
    (3,-1)*++!U\hbox{$R_7$},
    (2,-1)*{\bullet},
    (1,-1);(2,-1)**@{-},
    (2,-1)*++!U\hbox{$S_7$},
    (1,-1)*{\bullet},
    (0,0);(1,-1)**@{-},
    (1,-1)*++!U\hbox{$R_1$},
\end{xy}
$$
\caption{\label{AX}}
\end{figure}


Let $\mathcal S$ be a reducible plane sextic, decomposable as a cubic $\mathcal C$ with a cusp, a line $\mathcal L$ intersecting the cubic only at the cusp and a conic $\mathcal Q$ tangent to the line, (see figure \ref{plane sextic2}).



\begin{figure}[h]
    $$
  \begin{xy}
   <0cm,0cm>;<1.2cm,0cm>:
    (4,0.05)*++!D\hbox{$\mathcal L$},
    (0,0);(4,0)**@{-},
    (2.7,1)*++!D\hbox{$\mathcal C$},
    (2,0);(3.1,1.5)**\crv{(2.5,0)&(3.4,0.4)},
    (2,0);(3.1,-1.5)**\crv{(2.5,0)&(3.4,-0.4)},
   (1,1.1)*++!D\hbox{$\mathcal Q$},
   (1,0.6)*\xycircle(.6,.6){},
  \end{xy}
  $$
\caption{\label{plane sextic2}}
\end{figure}

\begin{prop}\cite{D1}
The map $X_0 \stackrel{\tau} \dashrightarrow Y$ factors as
$$\xymatrix{X_0 \ar@{-->}[d]^{/ \iota}\ar[r]^{/ \theta}& \mathbb P^ \times \mathbb P^1 \ar@{->}[d]^{\pi_1}\\ Y
\ar@{->}[r]^{\pi}& \mathbb P^2}$$
where the map $\pi_1$ is the double cover of $\mathbb P^2$ branched along the conic $\mathcal Q$ and the map $X_0 \stackrel{/ \theta} \to \mathbb P^1 \times \mathbb P^1$ is the double cover of $\mathbb P^1 \times \mathbb P^1$ branched along  $\pi_1^*(\mathcal C+\mathcal L)$.
\end{prop}

\noindent The map $\pi$ is given by the divisor $$\Delta= 9(L-E_0)-5E_{12}-4(E_{13}+E_{46}-3E_{56}-2(E_{24}+E_{36}+E_{45})-(E_{25}+E_{34})$$ and it contracts the $e_i$'s. Conversely,
\begin{prop}\cite{Naruki}
The canonical resolution of the double cover of $\mathbb P^2$ branched along the sextic $\mathcal S$ is a jacobian Kummer surface.
\end{prop}

\section{K3 Surfaces with $T_{X} \simeq U(2)\oplus U \oplus \langle -2 \rangle$}
\noindent Building on the results of the previous section and using the same notation, we describe the geometry of the K3 surfaces, appearing in our list, that have their period lattices isomorphic to $U(2)\oplus U \oplus \langle -2 \rangle$.\\

\noindent In the proposition \ref{San1}, we defined the curves $a_1, \cdots a_5$, which were components of the singular fibers of the elliptic fibration $Y \stackrel{f} \to \mathbb P^1$. Set $a_6:=e_6, \quad a_7:=e_7, \quad a_8:=E_{34}.$

\begin{lemma}\label{San2}
The $a_i$'s form an even eight on $Y$. If $X \stackrel{\mu} \dashrightarrow Y$ is the associated rational map of degree two, then $T_X \simeq U(2) \oplus U\oplus \langle -2 \rangle$.
\end{lemma}
\begin{proof}
A direct computation shows that
$$\sum_{i=1}^8 a_i=2C_{13}+ 2E_{34}+2\cdot \{
8(L-E_{0})-(5E_{12}+4E_{46}+3E_{13}+3E_{56}+E_{36}+E_{25}+E_{45})\}\in 2S_Y.$$
Recall that the rational double cover $ X \stackrel{\mu} \dashrightarrow Y $ associated to the even eight, $a_1, \cdots, a_8$, factors as
$$\xymatrix{Z \ar[r]^{\epsilon}\ar[d]_{p} & X
\ar@{-->}[dl]^{\mu}\\
Y}$$
where $Z \stackrel{p} \to Y$ is a degree two map branched along the $a_i$'s and $Z \stackrel{\epsilon} \to X$ contracts the
 eight exceptional curves $p^{-1}(a_i)$.
Let $X \stackrel{i} \to X$ be the covering Nikulin involution. We will describe the proper inverse image under $\mu$ of the elliptic fibration
$f:Y \stackrel{|D|} \to \mathbb P^1$ introduced in the proposition \ref{San1}. Each of the curves $e_1, e_2, e_3, e_4, e_5, e_8$ meet the $a_i$'s at exactly two points, thus their pullback under $p$ are smooth rational curves, each meeting an exceptional curve at two points. After blowing down the exceptional curves, they become six disjoint nodal curves $f_1, f_2, f_3, f_4, f_5, f_8$. It follows from the properties of $\mu^*$ that the divisor $\mu^*D$ is a nonzero effective nef divisor with self-intersection $0$. Thus the linear system $|\mu^*D|$ is of the form $|mD_0|, m\geq 1$, where $D_0$ is a smooth curve of genus one. Clearly the $f_i$'s are singular fibers of the elliptic fibration defined by $|D_0|$, $X \stackrel{g} \to \mathbb P^1$. Moreover as $i (f_i)=f_i$, the involution $i$ must act trivially on $\mathbb P^1$ under the map $g$.

\textit{Claim:} The general member of the linear system $|mD_0|$ is a smooth curve of genus one.
\textit{Proof of the claim:} Note that, since $mD_0 \cdot \mu^*C_{14} = \mu^*D \cdot \mu^*C_{14}=2 D \cdot C_{14} =2$, we must have $m \leq 2$. If $m=2$, then the preimage of the general member of $|D|$ would have to split into two disjoint smooth curves of genus one, $\mu^{-1}(D)= D_0+ D_1$, where $D_0 \stackrel{lin}\sim D_1,$ and $i(D_0)=D_1$. The last assertion contradicts the assumption that $i$ acts trivially on $\mathbb P^1$.

\noindent We deduce that the preimage of the $I_1$ fiber does not split into
two $I_1$ fibers but instead it becomes a singular fiber of type $I_2$. We denote its two components by  $S$ and $i^*S$. Finally, the $I^*_5$ fiber pulls back to a $I^*_{10}$-fiber whose components will be denoted by $$\arraycolsep=2pt\begin{array}{ll}E'_{14}+ E'_{15} + &2(C'_0+E'_{16}+C'_{16}+E'_{26}+C'_{12}+E'_{23}\\&+ i^*C'_{12}+i^*E'_{26}+i^*C'_{16}+ i^*E'_{16}+ i^*C'_0) + i^*E'_{14}+ i^*E'_{15}.\end{array}$$ The section $C_{14}$ meets the $a_i$'s at two points so it pulls back to a two-section $C'_{14}$.
We conclude that $X$ contains the following diagram of $(-2)$-curves. 

\begin{figure}[h] $$ \begin{xy}
<0cm,0cm>;<1.2cm,0cm>:
    (0,0)*{\bullet},
    (0,0)*++!D\hbox{$E'_{23}$},
    (0,0),(1,1)**@{-},
    (1,1)*{\bullet},
    (1,1);(2,1)**@{-},
    (1,1)*++!D\hbox{$C'_{12}$},
    (2,1)*{\bullet},
    (2,1);(3,1)**@{-},
    (2,1)*++!D\hbox{$E'_{26}$},
    (3,1)*{\bullet},
    (3,1);(4,1)**@{-},
    (3,1)*++!D\hbox{$C'_{26}$},
    (4,1)*++!D\hbox{$E'_{16}$},
    (4,1)*{\bullet},
    (4,1);(5,1)**@{-},
    (5,1)*++!D\hbox{$C'_{0}$},
    (5,1)*{\bullet},
    (6,1);(5,1)**@{-},
    (6,2)*{\bullet},
    (6,2)*++!D\hbox{$E'_{15}$},
    (5,1);(6,2)**@{-},
    (6,1)*++!D\hbox{$E'_{14}$},
    (6,1)*{\bullet},
    (7,0)*{\bullet},
    (6,1);(7,0)**@{-};
    (7,0);(6,-1)**@{-};
    (7,0)*++!D\hbox{$C'_{14}$},
    (8,1)*{\bullet},
    (8,1)*++!D\hbox{$S$},
    (7,0);(8,1)**@{-},
    (8,-1)*{\bullet},
    (8,-1)*++!U\hbox{$i^*S$},
    (7,0);(8,-1)**@{-},
    (8,1);(8,-1)**@{=},
    (6,-1)*{\bullet},
    (6,-1);(5,-1)**@{-},
    (6,-1)*++!U\hbox{$i^*E'_{14}$},
    (5,-1)*{\bullet},
    (5,-1);(4,-1)**@{-},
    (6,-2)*{\bullet},
    (6,-2)*++!U\hbox{$i^*E'_{15}$},
    (5,-1);(6,-2)**@{-},
    (5,-1)*++!D\hbox{$i^*C'_{0}$},
    (4,-1)*{\bullet},
    (3,-1);(4,-1)**@{-},
    (4,-1)*++!U\hbox{$i^*E'_{16}$},
    (3,-1)*{\bullet},
    (2,-1);(3,-1)**@{-},
    (3,-1)*++!U\hbox{$i^*C'_{16}$},
    (2,-1)*{\bullet},
    (1,-1);(2,-1)**@{-},
    (2,-1)*++!U\hbox{$i^*E'_{26}$},
    (1,-1)*{\bullet},
    (0,0);(1,-1)**@{-},
    (1,-1)*++!U\hbox{$i^*C'_{12}$},
\end{xy}
$$
\caption{\label{diagram2}}
\end{figure}

Let $M$ be the lattice of rank seventeen, generated by all the $(-2)$ curves of figure \ref{diagram2}. A
direct computation shows that $$D_{M} \simeq \mathbb Z/ 2\mathbb Z\oplus
\mathbb Z/ 4\mathbb Z \oplus \mathbb Z/ 4\mathbb Z.$$

\textit{Claim:} The lattice $S_{X}$ is two-elementary, i.e. $D_{S_X} \simeq (\mathbb Z/2\mathbb Z)^{a}$ for some $a \ge 0.$

\noindent Assuming the claim and letting $e$ be the index $[S_X:M]$, we get that
$$|D_M|=2^5=e^2 \cdot 2^a \Rightarrow e^2=2^{5-a}
\Rightarrow a=1,3 \textrm { or } 5. $$ If $a=5$, then $e=1$ and $M=S_{X}$, which contradicts the assumption that $S_{X}$ is two-elementary. If $a=1$, then by the classification of two-elementary lattices (\cite {N1}), $T_X$ would be isomorphic to $U \oplus U \oplus \langle -2 \rangle$. We have mentionned in the previous section that a K3 surface with such a period lattice contains finitely many smooth rational curves forming the figure \ref{AX}, which does not contain a subdiagram similar to the one in figure \ref{diagram2}. We conclude that $a=3$, or equivalentlty that $T_{X} \simeq U \oplus U(2) \oplus \langle -2 \rangle.$

\textit{Proof of the claim:} In order to prove the claim, it is enough to show the existence of an involution $X \stackrel{\theta} \to X$ satisfying $\theta^*_{S_X}=\mr{id}_{S_X}$ and $\theta^*_{T_X}=-\mr{id}_{T_X}$ (see proposition \ref{2-elem}). Observe that the curves $C_{13}$ and $C_{15}$ pull back to a curve of genus two $W,$ and a curve of genus one $E$, respectively, both invariant under $i$.

\noindent The cycle $$C'_{14}+E'_{14} +C'_0+E'_{16}+C'_{16}+E'_{26}+C'_{12}+E'_{23}+ i^*C'_{12}+i^*E'_{26}+i^*C'_{16}+ i^*E'_{16}+ i^*C'_0+ i^*E'_{14}$$ is disjoint from the genus one curve $E$. Consequently it defines the $I_{14}$ singular fiber of the elliptic fibration given by the linear system $|E|$, which we will denote by $v: X \stackrel{|E|}\to \mathbb P^1$. The curves $S$, $i^*S$, $E'_{15}$ and $i^*E'_{15}$ are sections of $v$. Fix $i^*E'_{15}$ as the zero section. Let $X \stackrel{\sigma} \to X$ be the map induced by the canonical involution on $X_{\eta}$, the generic fiber of $v$, i.e. $\sigma_{X_{\eta}}(x)=-x$ and let $X \stackrel{t_{E'_{15}}} \to X$ be the translation by the section $E'_{15}$.
Set $\gamma:=\sigma \circ t_{E'_{15}}$. The action of $\gamma$ on the period lattice is given by $$\gamma^*_{T_X}=\sigma^*_{T_X} \circ (t_{E'_{15}})^*_{T_X}=-\mr {id}_{T_X} \circ \mr{id}_{T_X}=-\mr{id}_{T_X}.$$ Moreover its action on the N\'eron-S\'ev\'eri lattice satisfies $\gamma^*_{S_X}=i^*_{S_X}.$ Indeed, to show this, we note that the divisor $S+i^*S -E'_{15}$ is linearly equivalent to the divisor
$$E'_{14}+ 2(C'_0+E'_{16}+C'_{16}+E'_{26}+C'_{12}+E'_{23}+ i^*C'_{12}+i^*E'_{26}+i^*C'_{16}+ i^*E'_{16}+ i^*C'_0)+ i^*E'_{14}+ i^*E'_{15}$$
which belongs to the kernel of the map $\beta$ of the theorem \ref{MW}. In other words $-S + E'_{15} = i^*S$ in $MW_{f_E}(X)$, or equivalently $\gamma^*S=i^*S$ and $\gamma^*E'_{15}=i^*E'_{15}$. It is now clear that $\gamma^*_M=i^*_M$ and consenquently that $\gamma^*_{S_X}=i^*_{S_X}.$
Define the map $\theta=i \circ \gamma$. By construction $\theta$ satisfies $\theta^*_{T_X
\oplus S_{X}}=-\mr {id}_{T_X} \oplus \mr{id}_{S_{X}}$, which implies, by the proposition \ref{2-elem}, that $S_X$ is two-elementary.
\end{proof}

\begin{definition}
Define the lattice $E_8^{twist}(-1)$ to be the lattice $E_8(-1)$
with the following modification $$ b(e_i, e_7)= \left\{
\begin{array}{cc} -4 & \textrm{ if } e_i = e_7\\ 2 & \textrm{ if }
e_i=e_6  \\ 0 & \textrm{ otherwise } \end{array} \right\} $$ where
the label of the vertices of the dynkin diagram associated to
$E_8(-1)$ is as follows

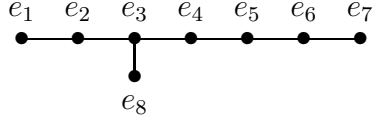
\begin{center}
$\begin{xy}
<0cm,0cm>;<0.5cm,0cm>:
(0,0)*{\bullet},
(0,0)*++!D\hbox{$e_1$},
(0,0),(1.5,0)**@{-}, (1.5,0)*{\bullet}, (1.5,0)*++!D\hbox{$e_2$}, (1.5,0);(3,0)**@{-},
(3,0)*{\bullet},(3,0)*++!D\hbox{$e_3$},(3,0);(4.5,0)**@{-}, (4.5,0)*{\bullet},(4.5,0)*++!D\hbox{$e_4$},(4.5,0);(6,0)**@{-},
(6,0)*{\bullet},(6,0)*++!D\hbox{$e_5$}, (6,0);(7.5,0)**@{-},(7.5,0)*{\bullet},(7.5,0)*++!D\hbox{$e_6$}, (7.5,0);(9,0)**@{-},
 (9,0)*{\bullet},(9,0)*++!D\hbox{$e_7$},
(3,-1)*{\bullet},(3,-1)*++!U\hbox{$e_8$},(3,0);(3,-1)**@{-},
\end{xy}$
\end{center}
Note that $|D_{E_8^{twist}(-1)}|=4.$
\end{definition}

\begin{prop}\label{twist} Let $Y$ be a general Kummer surface and $N$ be the Nikulin lattice (see definition \ref{Nikulin} in section 2) generated by the $a_i's$, then
there is a primitive embedding $E_8^{twist}(-1) \hookrightarrow N^{\perp}\subset S_Y.$
\end{prop}
\begin{proof}

The elements $E_{14}, C_{12}, E_{26}, C_{16}, E_{16}, C_0, E_{14}, 2C_{14}+a_5+a_8$ generate a sublattice of $N^{\perp}$ isomorphic to $E_8^{twist}(-1)$. To prove that the embedding is primitive, it is enough to show that this set can be completed to a basis of $N^{\perp}$. Let $T$ be the lattice generated by the elements $E_{14}, C_{12}, E_{26}, C_{16}, E_{16}, C_0, E_{14}, 2C_{14}+a_5+a_8$ and $2E_{23}+a_6+a_7$. $T$ is a sublattice of $N^{\perp}$ of finite index $e$. By explicitly taking the determinant of a matrix representing the bilinear form of $T$, we find that $|D_T|$=4 and therefore, $4=e^2| D_{N^{\perp}}|$. Either $e=1$ or $N^{\perp}$ is unimodular, but according to Milnor's classification \cite{Milnor}, there is no even unimodular lattice of rank nine. We conclude that $T=N^{\perp}$ and therefore the embedding is primitive.
\end{proof}

\begin{theorem}\label{San3}
Let $X$ be a $K3$ surface with $T_X \simeq U(2) \oplus U \oplus \langle -2 \rangle$. Then $X$ has a Nikulin involution such that if $X \stackrel{\mu} \dashrightarrow Y$ is the rational quotient map,

1) $Y$ is a Kummer surface,


2) there is a primitive embedding $E^{twist}_8(-1) \hookrightarrow  N^{\perp}.$

\end{theorem}

\noindent According to the theorem \ref{main thm}, we already know that $X$ admits a rational map of degree into a Kummer surface, so we need to show that the even eight associated to such a map is given by the $a_i$'s. Before proving this, we will show, using very geometrical arguments, that any K3 surface with $T_X \simeq U(2) \oplus U \oplus \langle -2 \rangle$ must contain a diagram of smooth rational curves meeting as in figure \ref{diagram2}. Let $\theta$ be the involution that satisfies $\theta^*_{T_{X}}=-\mr{id}_{T_{X}}$ and $\theta^*_{S_{X}}=\mr{id}_{S_{X}}$ which exists by proposition \ref{2-elem}. Nikulin showed that
$$X^{\theta}:= \{ x\in X\vert \theta(x)=x \}= E \textrm{ } \coprod_{i=1}^{7}\textrm{  }R_{i},$$ where $E$
is a curve of genus one and the $R_{i}$'s are disjoint rational
curves. The system $\vert E \vert$ defines an
elliptic fibration $X \stackrel{v} \to \mathbb P^{1}$, with a section $C$ and exactly one reducible singular fiber $\mathfrak F$ of type $I_{14}$. Moreover $\theta$ acts non trivially on $\mathbb P^1$ under $v$ (\cite{N4} p.1424-1430).

\noindent The $R_i$ are components of $\mathfrak F$. Since $\theta$ preserves all the rational curves, any section of $v$ must meet the singular fiber at a $R_i$ component. Thus the surface $X$ contains the following configuration of $(-2)$-curves

$$
  
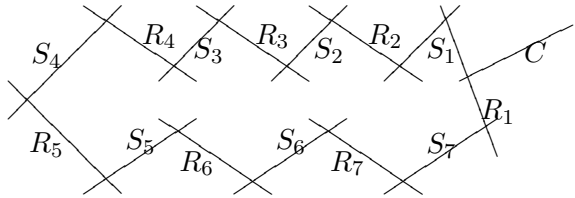
\begin{xy}
    <0cm,0cm>;<.5cm,0cm>:
    (0,0);(3,3)**@{-},
    (1,2.5)*++!U\hbox{\small$S_{4}$},
    (2,3);(5,1)**@{-},
    (4,3)*++!U\hbox{\small$R_{4}$},
    (4,1);(6,3)**@{-},
    (5.3,1)*++!D\hbox{\small$S_{3}$},
    (5,3);(8,1)**@{-},
    (7,3)*++!U\hbox{\small$R_{3}$},
    (7,1);(9,3)**@{-},
    (8.5,1)*++!D\hbox{\small$S_{2}$},
    (8,3);(11,1)**@{-},
    (10,3)*++!U\hbox{\small$R_{2}$},
    (10,1);(12,3)**@{-},
    (11.5,1)*++!D\hbox{\small$S_{1}$},
    (0,1);(3,-2)**@{-},
    (1,-1.5)*++!D\hbox{\small$R_{5}$},
    (2,-2);(5,0)**@{-},
    (3.5,0.2)*++!U\hbox{\small$S_{5}$},
    (4,0);(7,-2)**@{-},
    (5,-2)*++!D\hbox{\small$R_{6}$},
    (6,-2);(9,0)**@{-},
    (7.5,0.2)*++!U\hbox{\small$S_{6}$},
    (8,0);(11,-2)**@{-},
    (9,-2)*++!D\hbox{\small$R_{7}$},
    (10,-2);(13,0)**@{-},
    (11.5,0.1)*++!U\hbox{\small$S_{7}$},
    (11.5,3);(13,-1)**@{-},
    (13,1)*++!U\hbox{\small$R_{1}$},
    (12,1);(15,2.5)**@{-},
    (14,1)*++!D\hbox{\small$C$},
  \end{xy}
  $$

\noindent Consider the divisor $$F:=
2R_{2}+4S_{1}+6R_{1}+3C+5S_{7}+4R_{7}+3S_{6}+2R_{6}+S_{5}.$$
The system
$\vert F \vert$ defines another elliptic fibration, $X \stackrel{q}\to
\mathbb P^{1}$. By construction, $q$ has a singular fiber of type
$II^{*}$, a section $R_{5}$ and a 2-section $S_{2}$. Moreover, the curves $S_{4}, R_{4}, S_{3} \textrm { and } R_{3}$ must be
part of another reducible fiber of $q$, $\mathfrak J$. We will show by elimination that $\mathfrak J$ must be of type $I_{2}^{*}$. First, suppose that $\mathfrak J$ was of type $ I_{b} \textrm{ for } b \geq
4:$, then there would exist a rational smooth
curve $R$ meeting $S_{4}$ at exactly one point, $p$:

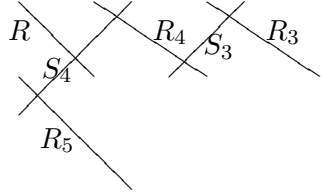
$$
\begin{xy}
    <0cm,0cm>;<.5cm,0cm>:
    (0,0);(3,3)**@{-},
    (1,2)*++!U\hbox{\small$S_{4}$},
    (2,1);(0,3)**@{-},
    (0,3)*++!U\hbox{\small$R$},
    (2,3);(5,1)**@{-},
    (4,3)*++!U\hbox{\small$R_{4}$},
    (4,1);(6,3)**@{-},
    (5.3,1)*++!D\hbox{\small$S_{3}$},
    (5,3);(8,1)**@{-},
    (7,3)*++!U\hbox{\small$R_{3}$},
    (0,1);(3,-2)**@{-},
    (1,-1.5)*++!D\hbox{\small$R_{5}$},
  \end{xy}
$$ Since $\theta^{*}$ acts as the identity on $S_{X}$,
$\theta(p)$ must also belong to $S_{4}\cap R$. But $R$ meets
$S_{4}$ at exactly one point, so $p$ must be fixed by $\theta$.
Either $p \in R_{5}\cap S_{4}\cap R$, which contradicts the fact
that $R_{5}$ is a section of $q$, or $p \in R_{4}\cap S_{4}\cap
R$, which contradicts the assumption that $S_{4}, R_{4}, S_{3}
\textrm { and } R_{3}$ are part of a fiber of type $I_{b}$. The fiber $\mathfrak J$ can neither be of type $I_{1}^{*}:$ Again, suppose it
is, then $S_{4}$ must be a component of multiplicity one ($S_{4}$ meets the section $R_{5}$). Denote by $R$ the other
component of multiplicity one of $\mathfrak J$ that meets $R_{4}$.
Two cases are possible:
\begin{enumerate}
\item $R$ does not meet the curve $S_{2}$:

$$
  
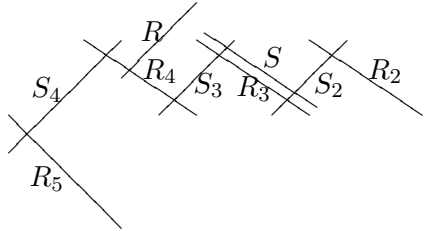
\begin{xy}
    <0cm,0cm>;<.5cm,0cm>:
    (0,0);(3,3)**@{-},
    (1,2.5)*++!U\hbox{\small$S_{4}$},
    (2,3);(5,1)**@{-},
    (4,3)*++!U\hbox{\small$R_{4}$},
    (3,2);(5,4 )**@{-},
    (3.8,4)*++!U\hbox{\small$R$},
    (4,1);(6,3)**@{-},
    (5.3,1)*++!D\hbox{\small$S_{3}$},
    (5,3);(8,1)**@{-},
    (6.5,.9)*++!D\hbox{\small$R_{3}$},
    (5.2,3.2);(8.2,1.2)**@{-},
    (7,3.3)*++!U\hbox{\small$S$},
    (7,1);(9,3)**@{-},
    (8.5,1)*++!D\hbox{\small$S_{2}$},
    (8,3);(11,1)**@{-},
    (10,3)*++!U\hbox{\small$R_{2}$},
    (0,1);(3,-2)**@{-},
    (1,-1.5)*++!D\hbox{\small$R_{5}$},
  \end{xy}
  $$

Then the sixth component of $\mathfrak J$, $S$ would
meet $S_{2}$ at exactly one point, $p$. Because $p \in S\cap S_{2}$, its image $\theta(p)$ belongs to $S\cap S_{2}$. So $\theta(p)=p$ and consequently $p \in R_{3}\cap
S_{2} \textrm { or } p \in R_{2}\cap S_{2}.$ The first case
contradicts the fact that $\mathfrak J$ is of type $I_{1}^{*}$, while in the second case, S would meet a component of
another singular fiber, which is absurd.
\item $R$ does meet the curve $S_{2}:$

$$
  
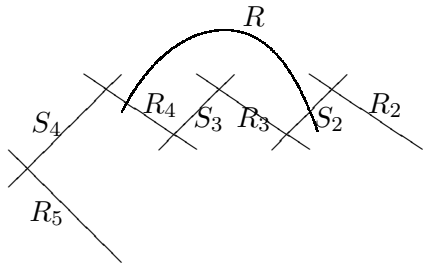
\begin{xy}
    <0cm,0cm>;<.5cm,0cm>:
    (0,0);(3,3)**@{-},
    (1,2.5)*++!U\hbox{\small$S_{4}$},
    (2,3);(5,1)**@{-},
    (4,3)*++!U\hbox{\small$R_{4}$},
    (3,2);(8.2,1.5)**\crv{(4,4)&(7,5)&(8,2)},
    (6.5,5.3)*++!U\hbox{\small$R$},
    (4,1);(6,3)**@{-},
    (5.3,1)*++!D\hbox{\small$S_{3}$},
    (5,3);(8,1)**@{-},
    (6.5,1)*++!D\hbox{\small$R_{3}$},
    (7,1);(9,3)**@{-},
    (8.5,1)*++!D\hbox{\small$S_{2}$},
    (8,3);(11,1)**@{-},
    (10,3)*++!U\hbox{\small$R_{2}$},
    (0,1);(3,-2)**@{-},
    (1,-1.5)*++!D\hbox{\small$R_{5}$},
  \end{xy}
  $$
Since $S_{2}$ is a 2-section of $q$ and $S_{2}\cdot R_{3}=1$,
necessarily $R \cdot S_{2}=1.$
 Let $p \in R\cap S_{2}$, then $\theta(p)=p$. It follows that either $p \in  R_{3}\cap
 S_{2},$ which contradicts again the assumption that $\mathfrak J$ is of type
$I_{1}^{*}.$  Or $p \in  R_{2}\cap S_{2}\cap R, $ which
implies again that $R$ meets a component of another singular fiber.
\end{enumerate}
We exclude also the case where the fiber $\mathfrak J$ is of type $I_{3}^{*}$ as follows. Suppose it is,
then $S_{4}$ has to be a component of multiplicity one and
$R_{4},S_{3}\textrm{ and }R_{3}$ components of multiplicity two.
Let $R$ be the fourth component of multiplicity two of $\mathfrak J$
and $L_{1}, L_{2}$ the two components of multiplicity one
intersecting it.

$$
  
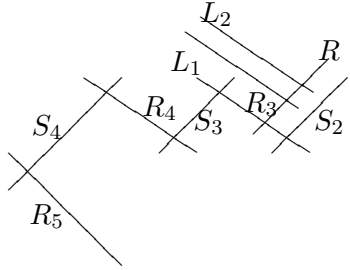
\begin{xy}
    <0cm,0cm>;<.5cm,0cm>:
    (0,0);(3,3)**@{-},
    (1,2.5)*++!U\hbox{\small$S_{4}$},
    (2,3);(5,1)**@{-},
    (4,3)*++!U\hbox{\small$R_{4}$},
    (4,1);(6,3)**@{-},
    (5.3,1)*++!D\hbox{\small$S_{3}$},
    (5,3);(8,1)**@{-},
    (6.7,3.1)*++!U\hbox{\small$R_{3}$},
    (6.5,1.5);(8.5,3.5)**@{-},
    (8.5,4.5)*++!U\hbox{\small$R$},
    (7.7,2.2);(4.7,4.2)**@{-},
    (5.5,5.5)*++!U\hbox{\small$L_{2}$},
    (8.1,2.6);(5.1,4.6)**@{-},
    (4.7,4.2)*++!U\hbox{\small$L_{1}$},
    (7,1);(9,3)**@{-},
    (8.5,1)*++!D\hbox{\small$S_{2}$},
    (0,1);(3,-2)**@{-},
    (1,-1.5)*++!D\hbox{\small$R_{5}$},
  \end{xy}
  $$

\noindent Since $L_{1}$ and $L_{2}$ do not meet any $R_{i}$'s nor any
$S_{i}$'s, they must be part of a reducible fiber of the fibration
defined by $\vert E \vert$. But this fibration is known to have an
unique reducible singular fiber. The fiber $\mathfrak J$ cannot be of type $I_{b}^{*}$, with $b \geq 4$, as $\mathfrak J$ would have $b+5$ components. But by theorem \ref{MW}, the number of components of $\mathfrak J$ cannot exceed 8, so $b \leq 3$. The same counting argument shows that
$\mathfrak J$ cannot be of type $II^{*}.$ The $\mathfrak J$ cannot be of type $III^{*}$. If it was, the group of sections of $q$ would be finite of
order $n$, where $n$ would satisfy by the theorem \ref{MW}
 $$2^{3}=|D_{S_{X}}|=\frac{1\cdot 2}{n^{2}}$$ but no integer satisfies this equation. Finally $\mathfrak J$ cannot be of type $IV^{*}$. Indeed in
this case one of the components with multiplicity one will have to
be disjoint from every $R_{i}$ and $S_{i}$. Such a component must
be part of a reducible fiber of the fibration defined by $\vert E
\vert$. But again, this fibration is known to have an unique
reducible singular fiber. By elimination, we have shown that $\mathfrak J$ must be of type
$I_{2}^{*}$. Moreover we have shown that $X$ contains the following
configuration of smooth rational curves

$$
  
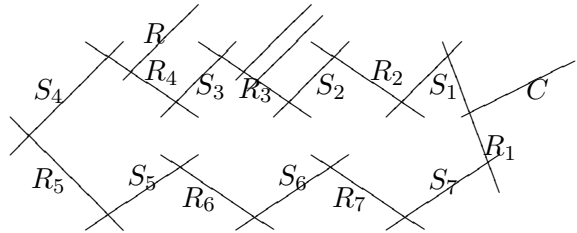
\begin{xy}
    <0cm,0cm>;<.5cm,0cm>:
    (0,0);(3,3)**@{-},
    (1,2.5)*++!U\hbox{\small$S_{4}$},
    (2,3);(5,1)**@{-},
    (4,3)*++!U\hbox{\small$R_{4}$},
    (3,2);(5,4 )**@{-},
    (3.8,4)*++!U\hbox{\small$R$},
    (4,1);(6,3)**@{-},
    (5.3,1)*++!D\hbox{\small$S_{3}$},
    (5,3);(8,1)**@{-},
    (6.5,.9)*++!D\hbox{\small$R_{3}$},
    (6,2);(8,4)**@{-},
    (6.3,1.7);(8.3,3.7)**@{-},
    (7,1);(9,3)**@{-},
    (8.5,1)*++!D\hbox{\small$S_{2}$},
    (8,3);(11,1)**@{-},
    (10,3)*++!U\hbox{\small$R_{2}$},
    (10,1);(12,3)**@{-},
    (11.5,1)*++!D\hbox{\small$S_{1}$},
    (0,1);(3,-2)**@{-},
    (1,-1.5)*++!D\hbox{\small$R_{5}$},
    (2,-2);(5,0)**@{-},
    (3.5,0.2)*++!U\hbox{\small$S_{5}$},
    (4,0);(7,-2)**@{-},
    (5,-2)*++!D\hbox{\small$R_{6}$},
    (6,-2);(9,0)**@{-},
    (7.5,0.2)*++!U\hbox{\small$S_{6}$},
    (8,0);(11,-2)**@{-},
    (9,-2)*++!D\hbox{\small$R_{7}$},
    (10,-2);(13,0)**@{-},
    (11.5,0.1)*++!U\hbox{\small$S_{7}$},
    (11.5,3);(13,-1)**@{-},
    (13,1)*++!U\hbox{\small$R_{1}$},
    (12,1);(15,2.5)**@{-},
    (14,1)*++!D\hbox{\small$C$},
  \end{xy}
  $$
By symmetry, the divisor class
$$F'=2R_{7}+2S_{7}+6R_{1}+3C+5S_{1}+4R_{2}+3S_{2}+2R_{3}+S_{3}$$
defines an elliptic fibration, $X \stackrel{h} \to \mathbb P^{1}$, for which the curves
$R_{6}+S_{5}+R_{5}+S_{4}$ are part of a reducible fiber of type $I_{2}^{*}$, $\mathfrak H:= S+S_{4}+2R_{5}+2S_{5}+2R_{6}+T+T'.$ Since $R \cdot
F'=0$, the curve $R$ is also part of a reducible fiber clearly
disjoint from $\mathfrak H$. It follows from theorem \ref{MW} that $R$ can only be
part of a fiber with two components so generically it is of type $I_{2}$. We conclude that the $R_i$'s, the $S_i$'s and the smooth rational curves $R, S, T$ define the diagram on $X$
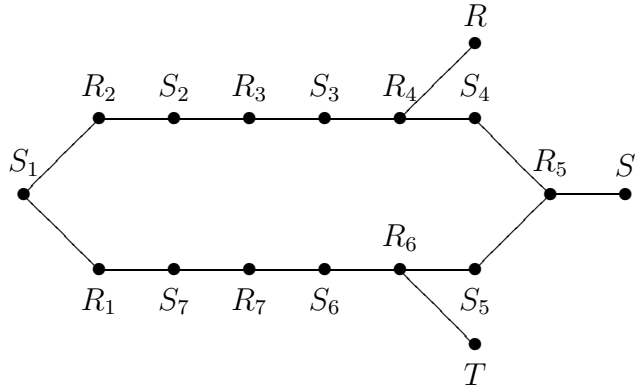
$$ \begin{xy}
<0cm,0cm>;<1cm,0cm>:
    (0,0)*{\bullet},
    (0,0)*++!D\hbox{$S_{1}$},
    (0,0),(1,1)**@{-},
    (1,1)*{\bullet},
    (1,1);(2,1)**@{-},
    (1,1)*++!D\hbox{$R_{2}$},
    (2,1)*{\bullet},
    (2,1);(3,1)**@{-},
    (2,1)*++!D\hbox{$S_{2}$},
    (3,1)*{\bullet},
    (3,1);(4,1)**@{-},
    (3,1)*++!D\hbox{$R_{3}$},
    (4,1)*++!D\hbox{$S_{3}$},
    (4,1)*{\bullet},
    (4,1);(5,1)**@{-},
    (5,1)*++!D\hbox{$R_{4}$},
    (5,1)*{\bullet},
    (6,1);(5,1)**@{-},
    (6,2)*{\bullet},
    (6,2)*++!D\hbox{$R$},
    (5,1);(6,2)**@{-},
    (6,1)*++!D\hbox{$S_{4}$},
    (6,1)*{\bullet},
    (7,0)*{\bullet},
    (6,1);(7,0)**@{-};
    (7,0);(6,-1)**@{-};
    (7,0)*++!D\hbox{$R_5$},
    (8,0)*{\bullet},
    (8,0)*++!D\hbox{$S$},
    (7,0);(8,0)**@{-},
    (6,-1)*{\bullet},
    (6,-1);(5,-1)**@{-},
    (6,-1)*++!U\hbox{$S_{5}$},
    (5,-1)*{\bullet},
    (5,-1);(4,-1)**@{-},
    (6,-2)*{\bullet},
    (6,-2)*++!U\hbox{$T$},
    (5,-1);(6,-2)**@{-},
    (5,-1)*++!D\hbox{$R_{6}$},
    (4,-1)*{\bullet},
    (3,-1);(4,-1)**@{-},
    (4,-1)*++!U\hbox{$S_6$},
    (3,-1)*{\bullet},
    (2,-1);(3,-1)**@{-},
    (3,-1)*++!U\hbox{$R_7$},
    (2,-1)*{\bullet},
    (1,-1);(2,-1)**@{-},
    (2,-1)*++!U\hbox{$S_7$},
    (1,-1)*{\bullet},
    (0,0);(1,-1)**@{-},
    (1,-1)*++!U\hbox{$R_1$},
\end{xy}
$$
The $I_{10}^*$ configuration apparent in this diagram defines an elliptic fibration for which $S$ is part of a $I_2$ fiber. By adding the other component of this fiber, we obtain a diagram of smooth rational curves meeting exactly as in the figure \ref{diagram2}.

We can now prove theorem \ref{San3}.

\begin{proof}
In order to explicitly construct the Nikulin involution $i$ we apply the inverse process used in lemma \ref{San2} to construct the involution $\theta$.  The curves $R$, $S$ and $T$, constructed above, are sections of the elliptic fibration $X \stackrel{v} \to \mathbb P^1$. Fix $T$ as the zero section. Let $\sigma$ be the map induced by inversion on the generic fiber and $t_{R}$ be the translation by $R$.  The map $\gamma = \sigma \circ t_R$ is an involution on $X$ such that $\gamma^*_{T_X}=-\mr{id}_{T_X}$.
Set $i := \theta \circ \gamma$.
As the maps $\theta^*$ and $\gamma^*$ commute on $\mr H^{2}(X,\mathbb{Z})$, $i$ is an involution. Moreover it follows from the fact that $i^{*}_{T_{X}}=\mr{id}_{T_{X}}$ that $i$ is a Nikulin involution.
Its eight fixed points can be easily described. Indeed a point $p$ is fixed by $i$ if and only if
$\gamma(p)=\theta (p)$. Since $\gamma$ acts trivially on the fibers of $v$ and $\theta$ acts non trivially on the same fibers, the fixed points of $i$ must lie on the only two fibers fixed by $\theta$, i.e $E$ and $\mathfrak F$. The points lying on $E$ satisfy $\gamma (p)=\theta(p)=p$ and therefore correspond to the four fixed point of $\gamma$ restricted to $E$.
The four remaining fixed points must lie on $\mathfrak F$. More precisely, two of them lie on $S_1$ and the other two lie on $R_5$. Indeed, since $S_{j} \cdot i^{*}(S_{j})=0$ for $j \neq 1$, a point $p$, lying on some $S_{j}'$s is fixed by $i$ if and only if $p\in S_{1}$ and $\theta(p)=p$. Thus the two fixed points of $\theta$ restricted to $S_1$ are fixed by $i$. Similarly, $R_{j} \cdot i^{*}(R_{j})=0$ for $j \neq 5$ implies that the two fixed points of $\gamma$ restricted to $R_5$ are fixed by $i$. Consider $X \stackrel{\mu} \dashrightarrow Y$ the rational quotient map associated to $i$ and denote by $A_1, \dots, A_8$ the induced even eight. The $A_i$'s and the image of the curves $S_i$, $R_i$, $S$, $R$ and $T$ under $\mu$ form the following diagram of rational curves on $Y$
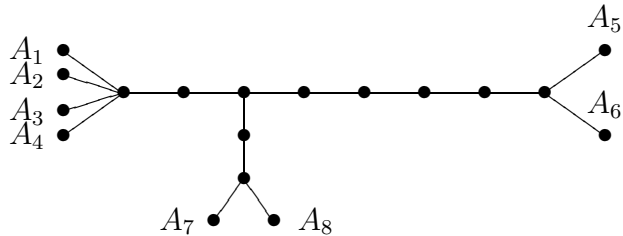
$$
  \begin{xy}
<0cm,0cm>;<0.8cm,0cm>:
    (0,0)*{\bullet},
    (0,0),(1,0)**@{-},
    (1,0)*{\bullet},
    (1,0);(2,0)**@{-},
    (2,0)*{\bullet},
    (2,0);(3,0)**@{-},
    (3,0)*{\bullet},
    (3,0);(4,0)**@{-},
    (4,0)*{\bullet},
    (4,0);(5,0)**@{-},
    (5,0)*{\bullet},
    (6,0);(5,0)**@{-},
    (6,0)*{\bullet},
    (7,0)*{\bullet},
    (7,0);(6,0)**@{-};
    (8,.71)*{\bullet},
    (7,0);(8,.71)**@{-},
    (8,.71)*++!D\hbox{$A_5$},
    (8,-.71)*{\bullet},
    (7,0);(8,-.71)**@{-},
    (8,-.71)*++!D\hbox{$A_6$},
    (-1,.71)*{\bullet},
    (-1,.71);(0,0)**@{-},
 (-1,.71)*++!R\hbox{$A_1$},
    (-1,.3)*{\bullet},
    (-1,.3);(0,0)**@{-},
    (-1,.3)*++!R\hbox{$A_2$},
    (-1,-.3)*{\bullet},
    (-1,-.3);(0,0)**@{-},
    (-1,-.3)*++!R\hbox{$A_3$},
    (-1,-.71)*{\bullet},
    (-1,-.71);(0,0)**@{-},
    (-1,-.71)*++!R\hbox{$A_4$},
    (2,-.71)*{\bullet},
    (2,0);(2,-.71)**@{-},
    (2,-1.42)*{\bullet},
    (2,-.71);(2,-1.42)**@{-},
    (1.5,-2.13)*{\bullet},
    (2,-1.42);(2.5,-2.13)**@{-},
    (1.5,-2.13)*++!R\hbox{$A_7$},
    (2.5,-2.13)*{\bullet},
    (2,-1.42);(1.5,-2.13)**@{-},
    (3.8,-2.13)*++!R\hbox{$A_8$}
\end{xy}
$$
The surface $Y$ admits therefore an elliptic fibration with a section and whose singular fibers are of the type $6I_2 + I^*_5 + I_1$. In the table of \cite{Sh1}, we find that the Mordell Weil group of such a fibration must be of order two. It follows now from the proposition \ref{San1} that $S_Y$ admits a primitive embedding of the \textit{Kummer lattice} which implies in turn by \cite{N2} that $Y$ is a Kummer surface. The primitive embedding $E_8^{twist}(-1) \hookrightarrow N^{\perp}$ follows from the proposition \ref{twist}.
\end{proof}

Let $\Sigma$ be a reducible plane sextic, decomposable as an irreducible quartic $\Gamma$ with an ordinary node and a cusp of type $A_4$,  a line $\Lambda_1$ passing trough the two singular points of $\Gamma$ and another line $\Lambda_2$ (see figure \ref{plane sextic3}).
\begin{figure}[h]
  $$
  \begin{xy}
    <0cm,0cm>;<1.2cm,0cm>:
    (0,-.3)*++!D\hbox{$\Lambda_2$},
    (4,-1)*++!D\hbox{$\Gamma$},
    (-.3,-1.15)*++!D\hbox{$\Lambda_1$},
     (0,-1);(1,1)**\crv{(1,-1)&(1.5,0)& (1,1)},
    (1,1);(3,0)**\crv{(0.5,-.5)&(2,-2)&(3,0)},
    (3,0);(4,-1)**\crv{(3,-1)&(4,-1)},
    (-.8,-.5);(4,1)**@{-},
    (-0.9,-1.3);(4,0.3)**@{-},
  \end{xy}
  $$
\caption{\label{plane sextic3}}
\end{figure}

\begin{theorem}\label{San4}
The rational quotient map $X \stackrel{\mu} \dashrightarrow Y$ decomposes as
$$\xymatrix{X \ar@{-->}[d]^{\mu}\ar[r]^{q_2}& Z \ar@{->}[d]^{q_1}\\ Y
\ar@{->}[r]^{q_0}& \mathbb P^2}$$ where $q_0$ is the canonical resolution of the double cover of $\mathbb P^2$ branched along $\Sigma$. The maps $q_1$ and $ q_2$ are the canonical resolutions of the double covers branched along $\Gamma$  and $q_1^*(\Lambda_1+\Lambda_2)$ respectively.
\end{theorem}
\begin{proof}
We have shown in the course of proving theorem \ref{San3} that the
map $\mu$ can be thought as the rational double cover associated
to the even eight defined by the $a_i$'s. Therefore there exists
on $X$ a curve of genus two defined by $W=\mu^*C_{13}$ that
satisfies  $$W \cap E= \{p_1, p_2, p_3, p_4 \}, \quad W \cap R_5 =
\{p_5, p_6\} \quad i(p_i)=p_i \quad i=1, \dots 6,$$ where the
$p_i$'s are fixed points of the Nikulin involution.  The quotient
$X/\theta$ is a smooth rational surface. Applying Hurwitz formula
we find that $K_{X/\theta}^2=-7$. Note that the image in $X/
\theta$ of the nine disjoint curves $S_i,$  for $i=1 \cdots 7$,
$R$ and $T$ have self intersection $-1$. We may therefore blow
them down to a surface $Z$, whose canonical divisor has
self-intersection 2. The map $X \stackrel{q_2} \to Z$ is ramified
along the $R_i$'s and $E$ and it sends $W$ to a smooth rational
curve $C$, with $C^2=2$. Since the involution $i$ commutes with
$\theta$ on $X$, it induces an involution $\bar{i}$ on $Z$. Note
that $\bar{i}$ has four fixed points on $C$, denoted by
$q_2(p_1),$ $q_2(p_2),$ $q_2(p_3),$ $q_2(p_4)$, thus $\bar{i}$
restricted to $C$ is the identity and the quotient map, $Z\to Z/
\bar{i}$, is ramified along $C$. Using Hurwitz formula a second
time, we find that $K_{Z/ \bar{i}}^2=6$. Moreover, under the
composition, $X \stackrel{q_2} \to Z \stackrel{/\bar{i}}\to Z/
\bar{i}$, the curves $R_2$, $R_3$ and $R_{5}$ map to rational
curves $R'_2, R'_3, R'_5$, with the intersection properties
$$R^{'2}_{5}=R^{'2}_2=-1, \quad R^{'2}_3=-2, \quad R'_5 \cdot
R'_2= R'_5 \cdot R'_3=0, \quad R'_2 \cdot R'_3=1.$$ We may
successively blow down $R'_5$, $R'_2$ and $R'_3$ to obtain a
rational surface whose canonical divisor has self intersection
equal to $9$, i.e. the projective plane. It is an easy exercise to
verify that under the composition, $X \stackrel{q_2} \to Z
\stackrel{q_1} \to \mathbb P^2$, the genus two curve $W$ maps to a
quartic $\Gamma$ with an ordinary node and a cusp of type $A_4$,
the $R_i$'s get contracted to a line $\Lambda_1$ passing through
the singularities of $\Gamma$, and the genus one curve $E$ maps to
another line $\Lambda_2$. The following proposition completes the
proof.
\end{proof}

\begin{prop}\label{San5}
A surface $Y$ is a general Kummer surface if and only if  $Y$ is the canonical resolution of a double cover of $\mathbb P^2$ branched along the reducible sextic $\Sigma = \Gamma + \Lambda_1 + \Lambda_2$.
\end{prop}
\begin{proof}
Assume first that $Y$ is a general Kummer surface. Recall that $Y$ admits a double plane model $Y \stackrel{\phi} \to \mathbb P^2$ that factors as $Y \stackrel{\tilde{\phi}} \to \tilde{\mathbb P}^2\stackrel{\epsilon}\to \mathbb P^2$ where $\tilde{\mathbb P}^2$ is the blowup of $\mathbb P^2$ at the fifteen singular points of the union of six lines (see figure \ref{plane sextic}).
The image in $\tilde{\mathbb P^2}$ of the curves $a_1, \cdots a_8, C_{15}, E_{15}, C_0, E_{16}, E_{26}, C_{12}, E_{26},$ $E_{14}$ $C_{14}$ and $C_{13}$ intersect as in the figure \ref{graph}. The numbers in the parenthesis denote the self-intersection.
\begin{figure}[h]
$$
  \begin{xy}
<0cm,0cm>;<1.25cm,0cm>:
    (0,0)*{\bullet},
    (0.1,0)*++!D\hbox{$\tilde C_{14}$\small{(-4)}},
    (0,0),(1,0)**@{-},
    (1,0)*{\bullet},
    (1,0);(2,0)**@{-},
    (1.1,0)*++!D\hbox{$\tilde E_{14}$\small{(-1)}},
    (2,0)*{\bullet},
    (2,0);(3,0)**@{-},
    (2,0)*++!D\hbox{$\tilde C_0$\small{(-4)}},
    (3,0)*{\bullet},
    (3,0);(4,0)**@{-},
    (3,0)*++!D\hbox{$\tilde E_{16}$\small{(-1)}},
    (4,0)*++!D\hbox{$C'_{16}$\small{(-4)}},
    (4,0)*{\bullet},
    (4,0);(5,0)**@{-},
    (5,0)*++!D\hbox{$\tilde E_{26}$\small{(-1)}},
    (5,0)*{\bullet},
    (6,0);(5,0)**@{-},
    (6,0)*++!D\hbox{$\tilde C_{12}$\small{(-4)}},
    (6,0)*{\bullet},
    (7,0)*{\bullet},
    (7,0);(6,0)**@{-};
    (7,0)*++!D\hbox{$\tilde E_{23}$\small{(-1)}},
    (8,0)*{\bullet},
    (7,0);(8,0)**@{-},
    (8,0)*++!D\hbox{$\tilde C_{23}$\small{(-2)}},
(0,0);(0,-1)**@{-},
(0,0);(-0.8,-1)**@{-},
(0,-1)*{\bullet},
 (-0.8,-1)*{\bullet},
(0,-1)*++!L\hbox{$\tilde a_{5}$\small{(-1)}},
(-0.8,-1)*++!R\hbox{$\tilde a_{8}$\small{(-1)}},
(2,-.71)*{\bullet},
    (2,0);(2,-.71)**@{-},
    (2,-.71)*++!L\hbox{$\tilde E_{15}$\small{(-1)}},
    (2,-1.42)*{\bullet},
    (2,-.71);(2,-1.42)**@{-},
    (2,-1.3)*++!L\hbox{$\tilde C_{15}$\small{(-4)}  },
    (2,-1.42);(0.7,-2.3)**@{-},
    (2,-1.42);(2.3,-2.3)**@{-},
    (2,-1.42);(3.3,-2.3)**@{-},
    (2,-1.42);(1.7,-2.3)**@{-},
    (2,-3.3);(0.7,-2.3)**@{-},
    (2,-3.3);(2.3,-2.3)**@{-},
    (2,-3.3);(3.3,-2.3)**@{-},
    (2,-3.3);(1.7,-2.3)**@{-},
(.7,-2.3)*{\bullet},
   (0.7,-2.3)*++!R\hbox{$\tilde a_{1}$\small{(-1)}},
(1.7,-2.3)*{\bullet},
    (1.85,-2.3)*++!R\hbox{$\tilde a_{2}$\small{(-1)}},
(2.3,-2.3)*{\bullet},
    (2.3,-2.3)*++!L\hbox{$\tilde a_{3}$\small{(-1)}},
(3.3,-2.3)*{\bullet},
    (3.3,-2.3)*++!L\hbox{$\tilde a_{4}$\small{(-1)}},
    (2,-3.3)*{\bullet},
(2,-3.3)*++!U\hbox{$\tilde C_{13}$\small{(-1)}},
(2,-3.3);(-0.5,-3.3)**@{-},
(-0.8,-1);(-0.5,-3.3)**@{-},
(0,-1);(-0.5,-3.3)**@{-},
(2,-3.3);(7,-3.3)**@{-},
(7,0);(7,-3.3)**@{-},
\end{xy}
$$
\caption{\label{graph}}
\end{figure}

\noindent The curve $\tilde C_{23}$ corresponds to the images $\tilde{\phi(a_6)}=\tilde{\phi (a_7)}$.
We may blow down the curves $\tilde a_{1}, \tilde a_{2}, \tilde a_{3}, \tilde a_{4}, \tilde a_5,$ $\tilde a_8, \tilde E_{14}, \tilde E_{16}, \tilde E_{15},\tilde E_{26}, \tilde E_{23}, \tilde C_{14}, \tilde C_{23}, \tilde C_{12}$ and $\tilde C_{16}$. This sequence of blow-downs defines a Cremona transformation
$$
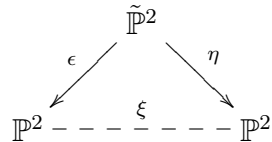
\xymatrix{ & \tilde{\mathbb P}^2 \ar[dl]_{\epsilon} \ar[dr]^{\eta} &
\\
\mathbb P^2 \ar@{--}[rr]^{\xi}& & \mathbb P^2}$$
which maps the six lines of figure \ref{plane sextic} as follows $$\xi(l_0)= \Lambda_1 \quad \xi(l_3)=\Gamma \quad \xi(l_5)=\Lambda_2$$
while the lines $l_2$, $l_4$ and $l_6$ get contracted to the singular points of $\Gamma$.

Conversely,
let $Y_0 \stackrel{\eta_1} \to \mathbb P^2$ be the blow-up of the $D_9$ singularity of $\Sigma$ (line passing through a cusp of type $A_4$) and denote by $\Upsilon_1$ the exceptional divisor. By its recursive construction (\cite{BPV} p.87), the canonical resolution of the double cover defined by $\Sigma$ is the canonical resolution of the double cover defined by $$\Sigma_1= \eta_1^* \Sigma - 2\Upsilon_1= \Gamma+ 2\Upsilon_1+ \Lambda_1 + \Upsilon_1 + \Lambda_2-2\Upsilon_1= \Gamma+ \Lambda_1 +\Lambda_2+\Upsilon_1$$ where we use the same symbol to denote a curve and its proper transform.
Recursively, if $Y_1 \stackrel{\eta_2} \to Y_0$ is the blow-up of the point $\Gamma \cap \Upsilon_1$, then the canonical resolution of the double cover defined by $\Sigma_1$ is the canonical resolution of the double cover defined by
 $$\Sigma_2= \eta_2^* \Sigma_1 - 2\Upsilon_2= \Gamma+ \Lambda_1 +\Lambda_2+\Upsilon_1 +\Upsilon_2.$$

\noindent Blow up the simple node of $\Gamma$ and denote the exceptional divisor by $\Upsilon_3$. We must now take the double cover branched along $$ \Sigma_3= \Gamma+ \Lambda_1 +\Lambda_2+\Upsilon_1 +\Upsilon_2 +\Upsilon_3.$$
Since the singularities of $\Sigma_3$ are now all of type $A_1$, the exceptional divisors coming from their resolution will not be part of the branch locus along which we will take the double cover. The resolution of $\Sigma_3$ gives a graph identical to the one of figure \ref{graph} in $\tilde{\mathbb P^2}$
where $\Gamma+ \Lambda_1 +\Lambda_2+\Upsilon_1 +\Upsilon_2 + \Upsilon_3$ corresponds exactly to $\tilde C_{13}+ \tilde C_{0}+ \tilde C_{15}+ \tilde C_{16}+ \tilde C_{12}+  \tilde C_{14}.$
\end{proof}



\section{K3 Surfaces with $T_X \simeq U(2)\oplus U(2) \oplus \langle -2 \rangle$}

The last section of this paper is very similar to the two previous ones. It aims to describe the geometry of the K3 surfaces with period lattice isomorphic to $U(2)\oplus U(2) \oplus \langle -2 \rangle$. 
\begin{prop}\label{San6}
Let $Y$ be a general Kummer surface. There exists an elliptic fibration on $Y$ whose singular fibers are of the
type $4I_2 + I^*_2 + I^*_1+I_1$ and whose Mordell-Weil group is
cyclic of order two.
\end{prop}
\begin{proof}
Consider the divisor
$$B=3(L-E_0)-2E_{12}-(E_{13}+E_{24}+E_{45}+E_{46}+E_{56}).$$
Geometrically, $B$ can be represented by the proper inverse image, under the map $Y \stackrel{\phi} \to \mathbb P^2$, of a cubic passing trough $p_{13}$, $p_{24}$, $p_{45}$, $p_{46}$, $p_{56}$ having a double point at $p_{12}$. By counting parameters, we immediately see that such a cubic exists. Generically, its proper transform is smooth and irreducible and satisfies $B^2=0$. Let $Y \stackrel{b} \to \mathbb P^1$ be the elliptic fibration defined by the linear system $|B|$.

The divisors
\begin{enumerate} \item $F_1=\underbrace{(L-E_0)-(E_{12}+E_{45})}_{b_1}+2(L-E_0)-(E_{12}+E_{13}+E_{24}+E_{46}+E_{56})$
\item $ F_2=\underbrace{2(L-E_0)-(E_{12}+E_{13}+E_{24}+E_{45}+E_{56})}_{b_2}+e_1$ \item
$F_3=\underbrace{3(L-E_0)-2E_{12}-(E_{13}+E_{24}+E_{36}+E_{45}+E_{46}+E_{56})}_{b_3}+E_{36}$\item
$F_4=\underbrace{E_{35}}_{b_4}+3(L-E_0)-2E_{12}-(E_{13}+E_{24}+E_{45}+E_{56}+E_{46}+E_{35})$
\end{enumerate}
define four fibers of type $I_2$. The divisor
$$F_5=\underbrace{(L-E_0)-(E_{12}+E_{56})}_{b_5}+E_{14}+2(C_0+E_{14}+C_{14})+E_{15}+E_{16}$$ defines a
fiber of type $I^*_2$. The divisor
$$F_6=e_6+e_7+2(C_{12}+E_{23})+E_{26}+E_{25}$$ defines a fiber of
type $I^*_1$.
Since $C_{15}\cdot B= C_{16}\cdot B =1$, we see that the rational curves $C_{15}$ and $C_{16}$ are sections of $b$. Moreover, as $\sum_{i=1}^6 \chi_{top}(F_i)=23$, there must be one more singular fiber $F_7$ of type $I_1.$ Finally, by theorem \ref{MW}, we deduce that $MW_{b}(X) \simeq \mathbb Z/2\mathbb Z.$
\end{proof}

Set $b_6:=e_6 \quad b_7:=e_7 \quad b_8:=E_{34}$.

\begin{lemma}\label{San7}
The $b_i$'s form an even eight on $Y$. If $X \stackrel{\nu} \dashrightarrow Y$ is the associated rational map of degree two, then $T_X \simeq U(2) \oplus U(2)\oplus \langle -2 \rangle$.
\end{lemma}

\begin{proof}
A direct computation shows that
$$\sum_{i=1}^8b_i=2C_{13}+2E_{35}+2E_{34}+2 \{
4(L-E_0)-(3E_{12}+2E_{45}+2E_{56}+E_{13}+E_{24}+E_{46})\}\in 2S_Y.$$

\noindent We will apply the same strategy as in lemma \ref{San2} and analyze the fibration, $b_{\nu}:X \stackrel{\nu^*B} \to \mathbb P^1$.
Let $X \stackrel{\iota} \to X $ be the associated Nikulin involution.
Likewise in lemma \ref{San2}, $\nu^*F_i$, for $i=1,2,3,4$, are nodal curves preserved by $\iota$. Hence $\iota$ preserves the fibers of $b_{\nu}$ ($\nu^*B$ cannot split). As for the remaining singular fibers, $\nu^*F_5$ becomes of type $I^*_4$, $\nu^*F_6$ becomes of type $I^*_2$ and $\nu^*F_7$ becomes of type $I_2.$
Finally $\nu^*C_{15}$ and $\nu^*C_{16}$ become 2-sections.. We conclude
that there exists an embedding of lattices of finite cokernel $U(2)\oplus D_8 \oplus D_6\oplus A_1 \hookrightarrow S_X$. Since $D_{U(2)\oplus D_8 \oplus D_6\oplus A_1}\cong (\mathbb Z/2\mathbb Z )^{\oplus 7}$, the overlattice $S_X$ must be a two-elementary lattice.
Note also that $b_\nu$ does not have a section. Indeed, if $C \cdot \nu^*D=1$, then $\iota^*C$ also would satisfy $$\iota^*C \cdot \iota^*(\nu^*D)=C \cdot \nu^* D=1.$$
Let $\tilde C \in S_Y$ be such that $\iota^*C+C=\nu^*\tilde C$. Then $2 \tilde C \cdot D= \nu^*\tilde C \cdot\nu^*D=2$, thus $\tilde C \cdot D =1.$ So either $C + \iota^*C= \nu^*C_{15}$ or $C + \iota^*C= \nu^*C_{16}.$
But we have seen that $C_{15}$ and $C_{16}$ do not split as they each meet the $b_i$'s at two points. This proves that $b_{\nu}$ does not have a section.
Let $J(X)$ be the jacobian fibration of $b_\nu$ (see \cite{CD} for the definition), then there is an embedding (\cite{Keum})$$T_{X} \hookrightarrow T_{J(X)}, \textrm{ with } T_{J(X)}/T_X \simeq \mathbb Z/2\mathbb Z.$$
Let $a$ be the length of $D_{S_X}$, i.e. $D_{S_X} \simeq (\mathbb Z/2\mathbb Z)^a.$ If $a=1$, then by the classification of two-elementary lattices (\cite{N1}), $T_X \simeq U \oplus U \oplus \langle -2 \rangle$ and $X$ contains finitely many rational curves meeting as in figure \ref{AX} which do not define any $I^*_4+I^*_2+I_2$ configuration. If $a=3$, then $T_{J(X)}$ has to be isomorphic to $U \oplus U \oplus \langle -2 \rangle$ and we get the same contradiction as in the previous case applied to the surface $J(X)$.
We conclude that $T_X \simeq U(2) \oplus U(2) \oplus \langle -2 \rangle.$
\end{proof}

\begin{theorem}\label{San8}
Let $X$ be a $K3$ surface with $T_X \simeq U(2) \oplus U(2) \oplus \langle -2 \rangle$. Then the even eight associated to the rational map $X \stackrel{\nu} \dashrightarrow Y$ is given by the $b_i$'s of lemma \ref{San7}.
\end{theorem}

\begin{proof}
From a purely lattice theoretical point of view, we deduce from lemma \ref{San7} that there exists an embedding $U(2)\oplus D_8 \oplus D_6\oplus A_1 \hookrightarrow S_X.$
Thus $X$ admits an elliptic fibration $X \stackrel{b_{\nu}} \to \mathbb P^1$ with singular fibers of type $I^*_4 + I^*_2 +I_2$, i.e. it contains the following diagram of smooth rational curves
$$
 
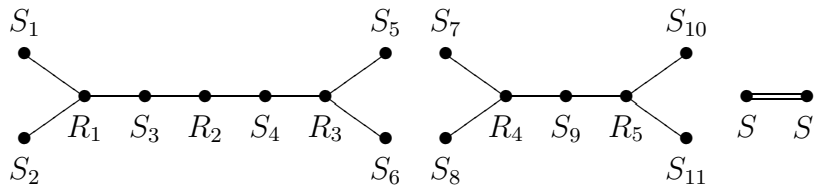
\begin{xy}
<0cm,0cm>;<0.8cm,0cm>:
    (0,0)*{\bullet},
    (0,0)*++!U\hbox{$R_1$},
    (0,0),(1,0)**@{-},
    (1,0)*{\bullet},
    (1,0)*++!U\hbox{$S_3$},
    (1,0);(2,0)**@{-},
    (2,0)*{\bullet},
    (2,0)*++!U\hbox{$R_2$},
    (2,0);(3,0)**@{-},
    (3,0)*{\bullet},
    (3,0);(4,0)**@{-},
    (3,0)*++!U\hbox{$S_4$},
    (4,0)*++!U\hbox{$R_3$},
    (4,0)*{\bullet},
    (6,0.71)*{\bullet},
(6,0.71)*++!D\hbox{$S_7$},
    (6,-0.71)*{\bullet},
(6,-0.71)*++!U\hbox{$S_8$},
    (10,0.71)*{\bullet},
(10,0.71)*++!D\hbox{$S_{10}$},
    (10,-0.71)*{\bullet},
(10,-0.71)*++!U\hbox{$S_{11}$},
    (6,0.71);(7,0)**@{-},
    (6,-0.71);(7,0)**@{-},
    (10,0.71);(9,0)**@{-},
    (10,-0.71);(9,0)**@{-},
    (7,0)*{\bullet},
    (8,0)*{\bullet},
(8,0)*++!U\hbox{$S_9$},
    (7,0);(8,0)**@{-};
    (7,0)*++!U\hbox{$R_4$},
    (9,0)*{\bullet},
 (9,0)*++!U\hbox{$R_5$},
    (9,0);(8,0)**@{-},
    (11,0)*{\bullet},
 (11,0)*++!U\hbox{$S$},
    (12,0)*{\bullet},
(12,0)*++!U\hbox{$S'$},
    (11,0);(12,0)**@{=},
    (-1,.71)*{\bullet},
 (-1,0.71)*++!D\hbox{$S_{1}$},
    (-1,.71);(0,0)**@{-},
    (-1,-.71)*{\bullet},
 (-1,-0.71)*++!U\hbox{$S_{2}$},
    (-1,-.71);(0,0)**@{-},
    (5,.71)*{\bullet},
 (5,0.71)*++!D\hbox{$S_{5}$},
    (5,.71);(4,0)**@{-},
    (5,-.71)*{\bullet},
 (5,-0.71)*++!U\hbox{$S_{6}$},
    (5,-.71);(4,0)**@{-},
    \end{xy}
$$
We shall see that the multisection index of $b_{\nu}$ is two. Note that by construction $b_{\nu}$ admits a $2$-section that we call $R_6$. Without loss of generality, $R_6 \cdot S_2= R_6 \cdot S_6=R_6 \cdot S_8= R_6 \cdot S_{11}=1.$
Let $X \stackrel{\theta} \to X$ be the involution acting trivially on $S_X$. Nikulin showed that $$X^{\theta} = \{ x\in X| \theta (x)=x\}$$ consists of the union of seven disjoint smooth rational curves \cite{N4}. In fact, the curves $R_1, R_2, R_3, R_4,$ $R_5, R_6$ belong to $X^{\theta}$. Indeed, the equalities $R_1\cdot S_1 = R_1\cdot S_2 =R_1 \cdot S_3 =1$ imply that there exist three points on $R_1$ fixed by $\theta$, consequently $R_1 \subset X^{\theta}$. The same argument applies for $R_3$, $R_4$ $R_5$ and $R_6$. Assume that $R_2 \notin X^{\theta}$ and denote by $p$ the point of $R_2 \cap S_3$. Because $R_2 \cdot S_3=1$, $\theta(p)=p$, so there exists $C \subset X^{\theta}$ such that $ p \in C \cap R_2 \cap S_3.$ If $F$ denotes the class of a fiber of $b_{\nu}$, then $C \cdot F \geq 4.$ The curves $S$ and $S'$ cannot be fixed by $\theta$ as they don't meet $S_1$, $S_5$, $S_7$ and $S_{10}$ which must intersect $X^{\theta}$. It implies that $C \cdot (S+S') \leq 4$, so $C \cdot (S+ S')=4$.
Similarly, if $q=R_2 \cap S_4$, then there exists $C' \subset X^{\theta}$, such that $C' \cdot F \geq 4.$
If $C \neq C'$, then $C' \cdot (S+S') =0$ which is absurd. If $C =C'$, then $C \cdot F \geq 8$, which contradicts $C \cdot (S + S') = 4$. So $R_2 \subset X^{\theta}.$ Denote by $R_7$ the remaining curve of $X^{\theta}$.
Necessarly $$(R_6+R_7)\cdot S= (R_6 + R_7)\cdot S'=2$$
and $R_7$ is also a $2$-section satisfying $$R_7 \cdot S_1= R_7 \cdot S_5=R_7 \cdot S_7= R_7 \cdot S_{10}=1.$$
Observe now that $$D_1= 2R_6+S_2+S_6+S_8+S_{11} \stackrel{lin}\sim 2R_7+S_1+S_5+S_7+S_{10}=D_2.$$ The divisors $D_1$ and $D_2$ are two $I^*_0$ fibers of an elliptic fibration $X \stackrel{h} \to \mathbb P^1.$ The curves $S_3, R_2, S_4$ and $S_9$ are components of the other singular fibers of $h$ and $R_1$, $R_3$, $R_4$, $R_5$ are sections of $h$. Fix $R_4$ as the zero section.
Let $\sigma$ be the map induced by inversion of the generic fiber and $t_{R_5}$ be the translation by $R_5$.  The map $\gamma = \sigma \circ t_5$ is an involution on $X$ such that $\gamma^*_{T_X}=-\mr{id}_{T_X}$.
Set $\iota := \theta \circ \gamma$.
As the maps $\theta^*$ and $\gamma^*$ commute on $\mr H^{2}(X,\mathbb{Z})$ and $\gamma^* \circ \theta^* _{T_X}=\mr{id}_{T_X}$ we conclude that $\iota$ is a Nikulin involution.
Note that $$R_1+S_3+R_2+S_4+R_3+S_3+R_7+S_1 \stackrel{lin}\sim R_4 + S_9+R_5+S_{11}+R_6+S_8$$ as they both define two disjoint $I_8$ and $I_6$ fibers of some other fibration. Therefore $-R_1 + R_5=R_3$ in $MW_h(X)$ and consequently $\iota^*R_1=\gamma^*R_1=R_3.$ Clearly $\iota^*R_4=\gamma^*R_4=R_5$. The involution $\gamma^*$ acts trivially on the curves $R_6$, $R_7$, $R_2$ and $S_9$. Consequently so does $\iota$ and each of those curves contains two of the eight fixed points of $\iota$. Finally, observe that $$\iota^*S_1=S_5, \quad \iota^*S_2=S_6, \quad \iota^*S_3=S_4, \quad \textrm{ and }\quad \iota^*S_7=S_{10} \quad \iota^*S_{8}=\iota^*S_{11}.$$
Consider $X \stackrel{\nu} \dashrightarrow Y$ the rational quotient map by $\iota$ and denote by $B_1, \dots B_8$ the induced even eight. Then $Y$ contains the following diagram of $-2$-curves

$$
 
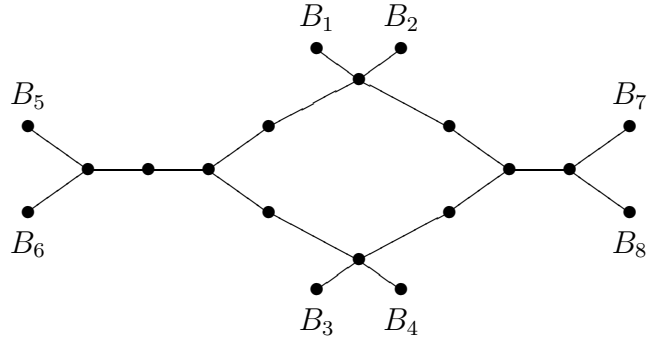
\begin{xy}
<0cm,0cm>;<0.8cm,0cm>:
    (0,0)*{\bullet},
    (0,0),(1,0)**@{-},
    (1,0)*{\bullet},
    (1,0);(2,0)**@{-},
   (2,0)*{\bullet},
   (6,0.71)*{\bullet},
    (6,-0.71)*{\bullet},
    (9,0.71)*{\bullet},
(9,0.71)*++!D\hbox{$B_7$},
    (9,-0.71)*{\bullet},
(9,-0.71)*++!U\hbox{$B_{8}$},
(4.5,1.5)*{\bullet},
(3.8,2)*{\bullet},
(5.2,2)*{\bullet},
(3.8,2)*++!D\hbox{$B_1$},
(5.2,2)*++!D\hbox{$B_2$},
(4.5,-1.5)*{\bullet},
(3.8,-2)*{\bullet},
(5.2,-2)*{\bullet},
(3.8,-2)*++!U\hbox{$B_3$},
(5.2,-2)*++!U\hbox{$B_4$},
 (4.5,1.5);(3,0.71)**@{-},
 (4.5,1.5);(6,0.71)**@{-},
(4.5,1.5);(3.8,2)**@{-},
(4.5,1.5);(5.2,2)**@{-},
(4.5,-1.5);(3,-0.71)**@{-},
 (4.5,-1.5);(6,-0.71)**@{-},
(4.5,-1.5);(3.8,-2)**@{-},
(4.5,-1.5);(5.2,-2)**@{-},
    (6,0.71);(7,0)**@{-},
    (6,-0.71);(7,0)**@{-},
    (9,0.71);(8,0)**@{-},
    (9,-0.71);(8,0)**@{-},
    (7,0)*{\bullet},
    (8,0)*{\bullet},
    (7,0);(8,0)**@{-};
    (-1,.71)*{\bullet},
 (-1,0.71)*++!D\hbox{$B_{5}$},
    (-1,.71);(0,0)**@{-},
    (-1,-.71)*{\bullet},
 (-1,-0.71)*++!U\hbox{$B_{6}$},
    (-1,-.71);(0,0)**@{-},
    (3,.71)*{\bullet},
    (3,.71);(2,0)**@{-},
    (3,-.71)*{\bullet},
    (3,-.71);(2,0)**@{-},
    \end{xy}
$$
By inspection of the diagram, we see that $Y$ admits an elliptic fibration with a section and whose singular fibers are of the type $4I_2 + I^*_2 + I^*_1$. In the table of \cite{Sh1}, we find that the Mordell Weil group of such a fibration must be of order two. It follows now from theorem \ref{San4} that $S_Y$ admits a primitive embedding of the \textit{Kummer lattice} which implies by Nikulin's criterion \cite{N2} that $Y$ is a Kummer surface. Finally, by proposition \ref{San6}, the $B_i$'s correspond to the $b_i$'s.
\end{proof}
\noindent Let $\mathcal X$ be a reducible plane sextic, decomposable as three lines $\mathcal L_1$, $\mathcal L_2$, $\mathcal L_3$ meeting at a point,  a line $\mathcal L_4$ and a conic $\mathcal W$ tangent to $\mathcal L_4$ (see figure \ref{plane sextic4}).
\begin{figure}[h]
  $$
  \begin{xy}
    <0cm,0cm>;<1cm,0cm>:
    (-1,0);(5.5,0)**@{-},
    (-1,-0.5);(5,3.5)**@{-},
    (4,3.5);(4,-0.5)**@{-},
    (4.3,3.4);(2.2,-0.5)**@{-},
    (1,0.9)*\xycircle(.9,.9){},
    (5.2,0.1)*++!\hbox{$\mathcal L_4$},
    (4.5,1)*++!\hbox{$\mathcal L_3$},
    (2.7,1)*++!\hbox{$\mathcal L_2$},
    (2.2,2)*++!\hbox{$\mathcal L_1$},
    (-.1,1)*++!\hbox{$\mathcal W$},
  \end{xy}
  $$
\caption{\label{plane sextic4}}
\end{figure}

\begin{lemma}\label{San9}
The map $X \stackrel{\nu} \dashrightarrow Y$ decomposes as
$$\xymatrix{X \ar@{-->}[d]^{\nu}\ar[r]^{\pi_2}& \mathbb P^1 \times \mathbb P^1 \ar@{->}[d]^{\pi_1}\\ Y
\ar@{->}[r]^{\pi}& \mathbb P^2}$$ where $\pi$ is the canonical resolution of the double cover of $\mathbb P^2$ branched along $\mathcal X$. The maps $\pi_1$ and $\pi_2$ are the canonical resolutions of the double covers branched along $\mathcal W$ and $\pi_1^*(\mathcal L_1 + \mathcal L_2 +\mathcal L_3 + \mathcal L_4)$ respectively.
\end{lemma}
\begin{proof}
Let $\mathbb P^1 \times \mathbb P^1  \stackrel{\pi_1}\to \mathbb P^2$ be the double cover branched along $\mathcal W$. Note that the line $\mathcal L_4$ splits under the cover and that $\pi_1^*(\mathcal L_1 + \mathcal L_2 + \mathcal L_3 + \mathcal L_4)$ is a divisor of bidegree $(4,4)$ and consequently it is an even divisor. The canonical resolution of the double cover defined by this even divisor is a K3 surface $X$ which admits a rational map $X \stackrel{\tilde \nu} \dashrightarrow Y$, where $Y$ is a general Kummer surface according to the proposition \ref{san} below. Applying the same techniques as in lemma \ref{San4}, we check that the corresponding even eight corresponds to the $b_i$'s of lemma \ref{San7} and therefore $\tilde \nu =\nu.$
 \end{proof}

\begin{prop}\label{san}
A surface $Y$ is a general Kummer surface if and only if  $Y$ is the canonical resolution of the double cover of $\mathbb P^2$ branched along the reducible sextic $\mathcal X$.
Moreover the double cover $Y \stackrel{\pi} \to \mathbb P^2$ contracts the $b_i'$s of lemma \ref{San7}.
\end{prop}

\begin{proof}
We proceed exactly as in proposition \ref{San5}. Consider the images of the curves $b_1, \cdots b_8,$ $C_{15}, E_{15}, C_0,$ $E_{14}, C_{14}, C_{13}, E_{25}, E_{16}, E_{23}, C_{12}, E_{26}, C_{16}$ under the map $Y \stackrel{\tilde {\phi}}\to \tilde{\mathbb P}^2$. Blow down the fourteen curves $\tilde b_1, \cdots, \tilde b_8,\tilde E_{14}, \tilde E_{15}, \tilde E_{23}, \tilde E_{26},\tilde E_{25},$ $\tilde E_{16}, \tilde C_{12}, \tilde C_{0}$. Since $\tilde{\mathbb P}^2$ is the blow-up of $\mathbb P^2$ at fifteen points, the resulting surface is, the blow-up of $\mathbb P^2$ at a point, i.e. the Hirzenbruch surface $\mathbb F_1$. Recall that $\mr{Pic}(\mathbb F_1)= \mathbb Z H \oplus \mathbb Z F$, where $H^2=1$, $H \cdot F=1$, and $F^2=0.$ Under the sequence of blow-downs, $\tilde C_{15}$ and $\tilde C_{16}$ become linearly equivalent to $H$ while $\tilde C_{14}$ and $\tilde C_{12}$ become linearly equivalent to $F$. Therefore the divisor $\tilde C_{14} - \tilde C_{15} \stackrel{lin} \sim E_p$, where is $E_p$ is the exceptional divisor of $\mathbb F_1$ (we use the same symbol to denote the curves $\tilde C_{13} $, $\tilde C_{14} $, $\tilde C_{15}$, $\tilde C_{16}$ and their images in $\mathbb F_1$). By keeping track of the intersections all along the blow-downs, we see that, on $\mathbb F_1$, $\tilde C_{13} \cdot \tilde C_{14}= \tilde C_{13} \cdot \tilde C_{15}=2$.  If follows that $\tilde C_{13} \cdot E_{p} =0$. Consider the contracting morphism $\mathbb F_1 \to \mathbb P^2$, under this map $\tilde C_{13}$ becomes a conic $\mathcal W$, the curves $\tilde C_{14}$, $\tilde C_{15}$, $\tilde C_{16}$ become three lines, $\mathcal L_1$, $\mathcal L_2$, $\mathcal L_3$, meeting at a point and the curve $\tilde C_{12}$ becomes a line $\mathcal L_4$ tangent to $\mathcal W$. Again we have constructed a Cremona transformation $$\xymatrix{ & \tilde{\mathbb P}^2 \ar[dl]_{\epsilon} \ar[dr]^{\eta} &
\\
\mathbb P^2 \ar@{--}[rr]^{\xi}& & \mathbb P^2}$$
which maps the six lines of the figure \ref{plane sextic} as follows $$\xi(l_3)=\mathcal W, \quad \xi(l_2)=\mathcal L_4, \quad \xi(l_4)=\mathcal L_1, \quad \xi(l_5)=\mathcal L_2, \quad \xi(l_6)=\mathcal L_3.$$
and the line $l_0$ gets contracted to the point of intersection $\mathcal L_1 \cap \mathcal L_2 \cap \mathcal L_3$.

Conversely, let $Y_0 \stackrel{\eta_1} \to \mathbb P^2$ be the blow-up of the point $p \in \mathcal L_1 \cap \mathcal L_2 \cap \mathcal L_3$ and denote by $\Upsilon_1$ the exceptional divisor. By construction, the canonical resolution of the double cover defined by $\mathcal X$ is the canonical resolution of the double cover defined by
$$ \mathcal X_1 = \eta^*\mathcal X - 2\Upsilon_1= \mathcal W + \mathcal L_1 + \mathcal L_2 + \mathcal L_3 + \mathcal L_4 + \Upsilon_1$$
(again we use the same symbol to denote a curve and its proper transform). Since the singularities of $\mathcal X_1$ are now all of type $A_1$, the exceptional divisors coming from their resolution will not be part of the branch locus along which we will take the double cover. Moreover, the canonical resolution of $\mathcal X_1$ gives rise to a graph of rational curves in $\tilde{\mathbb P}^2$ meeting exactly as the curves $ \tilde b_1, \cdots, \tilde b_8,\tilde E_{14}, \tilde E_{15}, \tilde E_{23}, \tilde E_{26},\tilde E_{25},$ $\tilde E_{16}, \tilde C_{12}, \tilde C_{0}, \tilde C_{13},\tilde C_{14}, \tilde C_{15}, \tilde C_{16}$. The divisor $\mathcal W + \mathcal L_1 + \mathcal L_2 + \mathcal L_3 + \mathcal L_4 + \Upsilon_1$ corresponds to $\tilde C_{13}+\tilde C_{14} +\tilde C_{15}+\tilde C_{16}+ \tilde C_{12} +\tilde C_{0}$.  It is clear now that the double cover along such a divisor is a general Kummer surface.
\end{proof}

\bibliographystyle{abbrv}
\bibliography{biblio}
\end{document}